\documentclass[12pt,reqno,twoside]{amsart}

\usepackage{amsfonts}
\usepackage{amsmath,amssymb,amsthm,amsthm}
\usepackage{mathtools}
\usepackage{dsfont}
\usepackage{etoolbox}
\usepackage{color}
\usepackage{bm}
\usepackage{subfigure}
\usepackage{units}
\usepackage{enumitem}

\usepackage{hyperref,graphicx}
\hypersetup{%
        colorlinks,
        breaklinks=true,
        plainpages=false,
        citecolor=.,
        linkcolor=.,
        urlcolor=.,
        bookmarksopen=true,%
        bookmarksnumbered=true,
        bookmarksdepth=8,
        unicode=true,%
        pdfpagelayout={SinglePage},pdffitwindow=true}%

\usepackage{cleveref}
\usepackage{booktabs}

\usepackage[dvips,bottom=1.4in,right=1in,top=1in, left=1in]{geometry}

\numberwithin{figure}{section}
\numberwithin{table}{section}

\makeatletter
\def\eqnarray{\stepcounter{equation}\let\@currentlabel=\theequation
\global\@eqnswtrue
\tabskip\@centering\let\\=\@eqncr
$$\halign to \displaywidth\bgroup\hfil\global\@eqcnt\z@
  $\displaystyle\tabskip\z@{##}$&\global\@eqcnt\@ne
  \hfil$\displaystyle{{}##{}}$\hfil
  &\global\@eqcnt\tw@ $\displaystyle{##}$\hfil
  \tabskip\@centering&\llap{##}\tabskip\z@\cr}

\def\endeqnarray{\@@eqncr\egroup
      \global\advance\c@equation\m@ne$$\global\@ignoretrue}

\usepackage{placeins}

\let\Oldsection\section
\renewcommand{\section}{\FloatBarrier\Oldsection}

 \usepackage[textsize=small]{todonotes}
\input{epsf}
\def\ms{\medskip}

\def\noi{\noindent}

\def\fvec{{\bf f}}

\def\svec{{\bf s}}

\def\uvec{{\bf u}}

\def\xvec{{\bf x}}

\def\utildevec{{\tilde {\bf u}}}

\def\Dmat{{\bf D}}

\def\Imat{{\bf I}}
\def\Jmat{{\bf J}}
\def\Kmat{{\bf K}}

\def\Mmat{{\bf M}}

\def\pmb#1{\setbox0=\hbox{$#1$}%
             \kern-.027em\copy0\kern-\wd0
             \kern+.009em\copy0\kern-\wd0
             \kern+.009em\copy0\kern-\wd0
             \kern+.009em\copy0\kern-\wd0
             \kern+.009em\copy0\kern-\wd0
             \kern+.009em\copy0\kern-\wd0
             \kern+.009em\copy0\kern-\wd0
             \kern-.045em\raise+.012em\copy0\kern-\wd0
             \kern+.009em\raise+.012em\copy0\kern-\wd0
             \kern+.009em\raise+.012em\copy0\kern-\wd0
             \kern+.009em\raise-.012em\copy0\kern-\wd0
             \kern+.009em\raise-.012em\copy0\kern-\wd0
             \kern-.018em\copy0\kern-\wd0\raise-.012em\box0}
\def\Pmb#1{\setbox0=\hbox{$#1$}%
             \kern-.033em\copy0\kern-\wd0
             \kern+.011em\copy0\kern-\wd0
             \kern+.011em\copy0\kern-\wd0
             \kern+.011em\copy0\kern-\wd0
             \kern+.011em\copy0\kern-\wd0
             \kern+.011em\copy0\kern-\wd0
             \kern+.011em\copy0\kern-\wd0
             \kern-.055em\raise+.015em\copy0\kern-\wd0
             \kern+.011em\raise+.015em\copy0\kern-\wd0
             \kern+.011em\raise+.015em\copy0\kern-\wd0
             \kern+.011em\raise-.015em\copy0\kern-\wd0
             \kern+.011em\raise-.015em\copy0\kern-\wd0
             \kern-.022em\copy0\kern-\wd0\raise-.015em\box0}
\def\alphavec{{\pmb \alpha}}

\vfill

\numberwithin{equation}{section}

\title[]{HIGH-FIDELITY DIGITAL TWINS: DETECTING AND LOCALIZING WEAKNESSES IN STRUCTURES}

\author[L\"ohner]{Rainald L\"ohner}
\address{Center for Computational Fluid Dynamics and 
Department of Physics and Astronomy, 
George Mason University, Fairfax, VA 22030, USA.}
\email{rlohner@gmu.edu}
\author[Airaudo]{Facundo N. Airaudo}
\address{Center for Computational Fluid Dynamics and 
Department of Physics and Astronomy, 
George Mason University, Fairfax, VA 22030, USA.}
\email{fairaudo@gmu.edu}
\author[Antil]{Harbir Antil}
\address{Center for Mathematics and Artificial Intelligence (CMAI) and 
Department of Mathematical Sciences, 
George Mason University, Fairfax, VA 22030, USA.}
\email{hantil@gmu.edu}
\author[W\"uchner]{Roland W\"uchner}
\address{Institut für Statik und Dynamik | Institute of Structural Analysis,
Beethovenstrasse 51, 38106 Braunschweig, Germany}
\email{r.wuechner@tu-braunschweig.de}
\author[Meister]{Fabian Meister}
\address{Institut für Statik und Dynamik | Institute of Structural Analysis,
Beethovenstrasse 51, 38106 Braunschweig, Germany}
\email{f.meister@tu-braunschweig.de}
\author[Warnakulasuriya]{Suneth Warnakulasuriya}
\address{Institut für Statik und Dynamik | Institute of Structural Analysis,
Beethovenstrasse 51, 38106 Braunschweig, Germany}
\email{s.warnakulasuriya@tu-braunschweig.de}

\begin{document}

\begin{abstract}
An adjoint-based procedure to determine weaknesses, or, more generally,
the material properties of structures is developed and tested. Given a
series of load cases and corresponding displacement/strain 
measurements, the material
properties are obtained by minimizing the weighted differences
between the measured and computed values. In a subsequent step,
techniques to minimize the number of load cases and sensors are
proposed and tested. \\
Several examples show the viability, accuracy and efficiency of
the proposed methodology and its potential use for
high fidelity digital twins.
\end{abstract}

\maketitle

\section{Introduction}
\label{sec:introduction}

Given that all materials exposed to the environment and/or undergoing 
loads eventually age and fail, the task of trying to detect and localize 
weaknesses in structures is common to many fields. To mention just a few:
airplanes, drones, turbines, launch pads and airport 
and marine infrastructure, bridges, high-rise buildings, wind turbines, 
and satellites.
Traditionally, manual inspection was the only way of carrying out this task,
aided by ultrasound, X-ray, or vibration analysis techniques.
The advent of accurate, abundant and cheap sensors, together with 
detailed, high-fidelity computational models
in an environment of digital twins has opened the possibility of 
enhancing and automating the detection and localization of 
weaknesses in structures. \\
From an abstract setting, it would seem that the task of determining
material properties from loads and measurements is an ill-posed
problem. After all, if we think of atoms, granules or some polygonal
(e.g. finite element [FEM]) subdivision of space, the amount of data
given resides in a much smaller space than the data sought.
If we think of a cuboid domain in $d$ dimensions with $N^d$ subdivisions,
the maximum amount of surface information/ data given is of 
$O(N^{d-1})$ while the data sought is of $O(N^d)$. \\
Another aspect that would seem to imply that this is an ill-posed
problem is the possibility that many different spatial distributions
of material properties could yield very similar or equal displacements
under loads. \\
On the other hand, the propagation of physical properties (e.g.
displacements, temperature, electrical currents, etc.) through the
domain obeys physical conservation laws, i.e. some partial differential
equations (PDEs). This implies that the material properties that can give
rise to the data measured on the boundary are restricted by these
conservation laws, i.e. are constrained. This would indicate that perhaps
the problem is not as ill-posed as initially thought. \\
As the task of damage detection is of such importance, many 
analytical techniques have been developed over the last decades 
\cite{cawley1979location,ladeveze1994modalfemupdate,maia1997localization,kim2004damage,rucka2006application,puelaubry2011meshadapt,mohan2013structural,chamoin2014updatingoffemmodels,mirzaee2015adjointaccel,alkayem2018structural,bunting2015eigenmodes,di2022data}.
Some of these were developed to identify weaknesses in
structures, others (e.g. \cite{ladeveze1994modalfemupdate}) to
correct or update finite element models.
The damage/weakness detection from measurements falls into
the more general class of inverse problems where material properties
are sought based on a desired cost functional
\cite{borrvall2003topology,lazarov2011filters,bunting2021novel,salloum2022optimization}. It is known that these inverse
problems are ill-defined and require regularization techniques. \\
The analytical methods depend on the measurement device
at hand, and one can classify broadly according to them.
The first class of analytical methods is based on changes
observed in (steady) displacements or strains
\cite{puelaubry2011meshadapt,chamoin2014updatingoffemmodels,alkayem2018structural,seidl2019forwbackw,di2022data}. The second
class considers velocities or accelerations in the time domain
\cite{kim2004damage,rucka2006application,mirzaee2015adjointaccel}.
The third class is based on changes observed in the frequency domain
\cite{cawley1979location,ladeveze1994modalfemupdate,maia1997localization,mohan2013structural,bunting2021novel}. \\
Some of the methods based on displacements, strains, velocities or
accelerations used adjoint formulations
\cite{troltzsch2010optimal,puelaubry2011meshadapt,chamoin2014updatingoffemmodels,mirzaee2015adjointaccel,antil2018frontiers,seidl2019forwbackw,lohner2020determination}
in order to obtain the gradient of the cost function with the least
amount of effort. In the present case, the procedures used are
also based on measured forces and displacements/strains, use adjoint
formulations and smoothing of gradients to quickly localize damaged
regions \cite{airaudo2023adjoint}. Unlike previous efforts, they are intended for
weakness/damage detection in the context of digital twins, i.e. we
assume a set of defined loadings
and sensors that accompany the structure (object, product, process)
throughout its lifetime in order to monitor its state. The digital
twins are assumed to contain finite element discretizations/models
of high fidelity, something that nowadays is common the aerospace 
industry. Therefore, the proposed approach fits well into the overall
workflow of high-level CAD environments and high fidelity FEM models
seen in the design phase.

\section{Assumptions}
\label{sec:assumptions}

What follows relies on the following set of assumptions:
\begin{itemize}
\item[-] Monitoring the weakening of a structure is carried out
by applying a set of $n$ different forces $\fvec_i, i=1,n$ and measuring
the resulting displacements $\uvec^{md}_{ij}, i=1,n, j=1,m$ and/or
strains $\svec^{ms}_{ij}, i=1,n, j=1,m$ at $m$ different
locations $\xvec_j, j=1,m$ (the intrinsic assumption is that the
forces can be standardized and perhaps even maintained throughout the 
life of the structure);
\item[-] The weakening of a structure may occur at any location, i.e.
there are no regions that are excluded for weakening; this is the
most conservative assumption, and could be relaxed under certain
conditions;
\item[-] The sensors for displacements and strains are limited in
their ability to measure by noise/signal ratios, i.e. actual
displacements and strains have to be larger than a certain threshold
to be of use:
\begin{equation} \label{eq:state_constraints}
    |\uvec^m| \ge u_0 ~~,~~ |\svec^m| \ge s_0 ~~.
\end{equation}
\item[-] The type of force used to monitor the weakening of a structure
is limited by practical considerations; this implies that the number
of different forces is limited, and can not assume arbitrary distributions
in space.
\item[-] The weakening of a structure may be described by a field
$\alpha(\xvec)$, where $0 \le \alpha(\xvec) \le 1$ and $\alpha(\xvec)=0$
corresponds to total failure (no load bearing capability) while 
$\alpha(\xvec)=1$ is the original state;
\item[-] The displacements, strains and stresses of the structure
are well
described by a sufficiently fine finite element discretization
(e.g. trusses, beams, plates, shells, solids) \cite{zienkiewicz2005finite,simo2006computational}, which
results in a system of equations for each load case:
\begin{equation} \label{eq:forward}
    \Kmat \uvec_i = \fvec_i ~~,~~ i=1,n
\end{equation}
\item[ ] where $\uvec_i$ are the displacements and $\Kmat$ the usual
stiffness matrix, which is obtained by assembling all the
element matrices:
\begin{equation} \label{eq:stiffness_assembly}
    \Kmat = \sum_{e=1}^{N_e} \Kmat_{e}  ~~.
\end{equation}
\end{itemize}

\section{Determining material properties via optimization}
\label{sec:determining_material_prop}

The determination of material properties (or weaknesses)
may be formulated as an 
optimization problem for the strength factor $\alpha(\xvec)$ as follows: 
Given $n$ force loadings $\fvec_i, i=1,n$ and $n \cdot m$ corresponding
measurements at
$m$ measuring points/locations $\xvec_j, j=1,m$ of their
respective displacements $\uvec^{md}_{ij}, i=1,n, j=1,m$ or strains
$\svec^{ms}_{ij}, i=1,n, j=1,m$, obtain the spatial distribution 
of the strength factor $\alpha$ that minimizes the cost function:
\begin{equation} \label{eq:objective}
    I(\uvec_{1,..,n},\svec_{1,..,n},\alpha) = 
  {1 \over 2} \sum_{i=1}^n \sum_{j=1}^m w^{md}_{ij} 
             ( \uvec^{md}_{ij} - \Imat^d_{ij} \uvec_i )^2
+ {1 \over 2} \sum_{i=1}^n \sum_{j=1}^m w^{ms}_{ij} 
             ( \svec^{ms}_{ij} - \Imat^s_{ij} \svec_i )^2
                                                   ~~,
\end{equation}
subject to the finite element description (e.g. trusses,
beams, plates, shells, solids) of the structure \cite{zienkiewicz2005finite,simo2006computational} under 
consideration (i.e. the digital twin/system \cite{maia1997localization,chinesta2020virtual}):
\begin{equation} \label{eq:forward2}
    \Kmat \cdot \uvec_i = \fvec_i ~~,~~ i=1,n
\end{equation}
where $w^{md}_{ij}, w^{ms}_{ij}$ are displacement and strain weights, 
$\Imat^d, \Imat^s$ interpolation matrices that are used to obtain 
the displacements and strains from the finite element mesh
at the measurement locations, and $\Kmat$ the usual stiffness matrix, 
which is obtained by assembling all the element matrices:
\begin{equation} \label{eq:strf_assembly}
    \Kmat = \sum_{e=1}^{N_e} \alpha_{e} \Kmat_{e}  ~~,
\end{equation}
where the strength factor $\alpha_{e}$ of the elements has already been
incorporated. We note in passing that in order to ensure that $\Kmat$
is invertible and non-degenerate $\alpha_{e} > \epsilon > 0$.
Note that the optimization problem given by Eqns.\eqref{eq:objective}-\eqref{eq:strf_assembly} does not 
assume any specific choice of finite element basis functions, i.e. 
is widely applicable.

\subsection{Optimization via adjoints} \label{subsec:optimization_via_adjoints}

The objective function can be extended to the Lagrangian functional
\begin{equation} \label{eq:lagrangian}
    L(\uvec_{1,..,n},\alpha,\utildevec_{1,..,n}) = I(\uvec_{1,..,n},\alpha)
+ \sum_{i=1}^n \utildevec^t_i ( \Kmat \uvec_i - \fvec_i ) ~~,
\end{equation}
where $\utildevec_i$ are the Lagrange multipliers (adjoints).
Variation of the Lagrangian with respect to each of the measurements
then results in:
\begin{subequations}
\begin{eqnarray} \label{eq:lagrangian_derivatives}
    {{dL} \over {d\utildevec_i}} = \Kmat \uvec_i - \fvec_i = 0 \\
    {{dL} \over {d\uvec_i}}      = 
    \sum_{j=1}^m w^{md}_{ij} \Imat^d_{ij} (\uvec^{md}_{ij} - \Imat^d_{ij} \uvec_i)
    + \sum_{j=1}^m w^{ms}_{ij} \Jmat^s_{ij} (\svec^{ms}_{ij} - \Imat^s_{ij} \svec_i)
    + \Kmat^t \utildevec_i = 0 \\
    {{dL} \over {d\alpha_e}}     = 
    \sum_{i=1}^n \utildevec_i^t {{d\Kmat} \over {d\alpha_e}} 
                         \uvec_i =
    \sum_{i=1}^n \utildevec_i^t \Kmat^e \uvec_i ~~,
\end{eqnarray}
\end{subequations}
where $\Jmat^s_{ij}$ denotes the relationship between the displacements
and strains (i.e. the derivatives of the displacement field on the
finite element mesh and the location $\xvec_j$ (see Section \ref{sec:interpolation} below)). \\
The consequences of this rearrangement are profound:
\begin{itemize}
\item[-] The gradient of $L$, $I$ with respect to $\alpha$ may be 
obtained by solving $n$ forward and adjoint problems; i.e.
\item[-] Unlike finite difference methods, which require at least $n$ 
forward problems per design variable, the number of forward and adjoint 
problems to be solved is {\bf independent of the number of variables 
used for $\alpha$} (!);
\item[-] Once the $n$ forward and adjoint problems have been solved, 
the cost for the evaluation of the gradient of each design variable 
$\alpha_e$ only involves the degrees of freedom of the element, i.e. 
is of complexity $O(1)$;
\item[-] For most structural problems $\Kmat=\Kmat^t$, so if a direct 
solver has been employed for the forward problem, the cost for the 
evaluation of the adjoint problems is negligible;
\item[-] For most structural problems $\Kmat=\Kmat^t$, so if an iterative 
solver is employed for the forward and adjoint problems, the 
preconditioner can be re-utilized.
\end{itemize}

\subsection{Optimization steps} \label{subsec:optimization_steps}

\ms \noi
An optimization cycle using the adjoint approach is composed 
of the following steps: \\
For each force/measurement pair $i$:
\begin{itemize}
\item[-] With current $\alpha$: solve for the displacements 
$\rightarrow \uvec_i$
\item[-] With current $\alpha$, $\uvec_i$ and 
$\uvec^{md}_{ij}, \svec^{md}_{ij}$: solve for the adjoints 
$\rightarrow \utildevec_i$
\item[-] With $\uvec_i, \utildevec_i$: obtain gradients
$\rightarrow I^i_{,\alpha}=L^i_{,\alpha}$
\end{itemize}
\par \noi
Once all the gradients have been obtained:
\begin{itemize}
\item[-] Sum up the gradients 
$\rightarrow I_{,\alpha} = \sum_{i=1}^n I^i_{,\alpha}$
\item[-] If necessary: smooth gradients $\rightarrow I^{smoo}_{,\alpha}$
\item[-] Update $\alpha_{new} = \alpha_{old} - \gamma I^{smoo}_{,\alpha}$.
\end{itemize}
\noi
Here $\gamma$ is a small stepsize that can be adjusted so as to obtain
optimal convergence (e.g. via a line search method).

\section{Interpolation of displacements and strains}
\label{sec:interpolation}

The location of a displacement or strain gauge may not coincide
with any of the nodes of the finite element mesh. Therefore, 
in general, the displacement $\uvec_k$ at a measurement location 
$\xvec^m_k$ needs to be obtained via the interpolation matrix 
$\Imat^d_k$ as follows:
\begin{equation} \label{eq:displ_interp}
    \uvec_k(\xvec^m_k) = \Imat^d_k(\xvec^m_k) \cdot \uvec ~~,
\end{equation}
where $\uvec$ are the values of the displacements vector at all 
gridpoints. \\
In many cases it is much simpler to install strain gauges instead of
displacement gauges. In this case, the strains need to be obtained
from the displacement field. This can be written formally as:
\begin{equation} \label{eq:strain}
    \svec = \Dmat \uvec ~~,
\end{equation}
where the `derivative matrix' contains the local values of the
derivatives of the shape-functions of $\uvec$.
The strain at an arbitrary position $\xvec^m_i$ is obtained via 
the interpolation matrix $\Imat^s_k$ as follows:
\begin{equation} \label{eq:strain_interp}
    \svec_k(\xvec^m_k) = \Imat^s_k(\xvec^m_k) \cdot \svec 
                      = \Imat^s_k(\xvec^m_k) \cdot \Dmat \cdot \uvec ~~.
\end{equation}

Note that in many cases the strains will only be defined in the 
elements, so that the interpolation matrices for displacements 
and strains may differ.

\section{Choice of weights}
\label{sec:weights}

The cost function is given by Eqn.\eqref{eq:objective}, repeated here for clarity:
\begin{equation} \label{eq:objective2}
    I(\uvec_n,\alpha) = 
  {1 \over 2} \sum_{i=1}^n \sum_{j=1}^m w^{md}_{ij} 
             ( \uvec^{md}_{ij} - \Imat^d_{ij} \uvec_i )^2
+ {1 \over 2} \sum_{i=1}^n \sum_{j=1}^m w^{ms}_{ij} 
             ( \svec^{ms}_{ij} - \Imat^s_{ij} \svec_i )^2
                                                   ~~.
\end{equation}

One can immediately see that the dimensions of displacements and
strains are different. This implies that the weights should be chosen
in order that all the dimensions coincide. The simplest way of
achieving this is by making the cost function dimensionless. This
implies that the displacement weights $w^{md}_{ij}$ should be of
dimension [1/(displacement*displacement)] (the strains are already
dimensionless). Furthermore, in order to make the procedures more
generally applicable they should not depend on a particular
choice of measurement units (metric, imperial, etc.). This implies
that the weights for the displacements and strains should be of
the order of the characteristic or measured magnitude.
Several options are possible, listed below.

\subsection{Local Weighting} \label{subsec:local_weighting}
In this case 
\begin{equation} \label{eq:local_weighting}
   w^{md}_{ij}={1 \over{(\uvec^{md}_{ij})^2}} ~~;~~ 
   w^{ms}_{ij}={1 \over{(\svec^{ms}_{ij})^2}} ~~; 
\end{equation}
this works well, but may lead to an `over-emphasis' of small 
displacements/strains that are in regions of marginal interest.

\subsection{Average weighting} \label{subsec:average_weighting}
In this case one first obtains the average
of the absolute value of the displacements/strains for a loadcase and uses
them for the weights, i.e.: 
\begin{equation} \label{eq:average_weighting}
    u_{av}={{\sum_{j=1}^m |\uvec^{md}_{ij}|} \over m} ~~; ~~
   w^{md}_{ij}={1 \over u^2_{av}} ~~;
   s_{av}={{\sum_{j=1}^m |\svec^{ms}_{ij}|} \over m} ~~; ~~
   w^{ms}_{ij}={1 \over s^2_{av}} ~~;
\end{equation}
this works well, but may lead to an `under-emphasis' of small 
displacements/strains that may occur in important regions.

\subsection{Max weighting} \label{subsec:max_weighting}
In this case one first obtains the maximum
of the absolute value of the displacements/strains for a loadcase and uses
them for the weights, i.e.:
\begin{subequations}
\begin{eqnarray} \label{eq:max_weighting}
    u_{max}=max(|\uvec^{md}_{ij}|, j=1,m) ~~;~~
   w^{md}_{ij}={1 \over u^2_{max}} ~~; \\
   s_{max}=max(|\svec^{ms}_{ij}|, j=1,m) ~~; ~~
  w^{ms}_{ij}={1 \over s^2_{max}} ~~;
\end{eqnarray}
\end{subequations}
this also works well for many cases, but may lead to an 
`under-emphasis' of smaller displacements/strains that can occur in 
important regions.

\subsection{Local/Max weighting} \label{subsec:local_max_weighting}
In this case
\begin{equation} \label{eq:local_max_weighting}
    w^{md}_{ij}={1 \over {max(\epsilon u_{max}, |\uvec^{md}_{ij}|))^2}} ~~;~~ 
   w^{ms}_{ij}={1 \over {max(\epsilon s_{max}, |\svec^{ms}_{ij}|))^2}} ~~;
\end{equation}
with $\epsilon=O(0.01-0.10)$; this seemed to work best of all, 
as it combines local weighting with a max-bound minimum for local values.

\section{Smoothing of gradients}
\label{sec:smoothing}

The gradients of the cost function with respect to $\alpha$ allow for
oscillatory solutions. One must therefore smooth or `regularize' the
spatial distribution. This happens naturally when using few degrees of
freedom, i.e. when $\alpha$ is defined via other spatial shape functions
(e.g. larger spatial regions of piecewise constant $\alpha$ \cite{puelaubry2011meshadapt}). 
As the (possibly oscillatory) gradients obtained in the (many) finite 
elements are averaged over spatial regions, an intrinsic smoothing occurs.
This is not the case if $\alpha$ and the gradient are defined and 
evaluated in each element separately, allowing for the largest degrees of
freedom in a mesh and hence the most accurate representation.
Three different types of smoothing or `regularization' were considered.
All of them start by performing a volume averaging from elements
to points:
\begin{equation} \label{eq:volume_avg}
    \alpha_p = {{ \sum_{e} \alpha_{e} V_{e} } \over
               { \sum_{e} V_{e} }} ~~,
\end{equation}
where $\alpha_p, \alpha_{e}, V_{e} $ denote the value of $\alpha$ at
point $p$, as well as the values of $\alpha$ in element $e$ and the
volume of element $e$, and the sum extends over all the elements
surrounding point $p$.

\subsection{Simple point/element/point averaging} \label{subsec:simple_avg}
In this case, the values of $\alpha$ are cycled between elements 
and points. When going from point values to element values, 
a simple average is taken:
\begin{equation} \label{eq:simple_avg}
    \alpha_{e} = { 1 \over n_{e} } \sum_i \alpha_i ~~,
\end{equation}
where $n_{e}$ denotes the number of nodes of
an element and the sum extends over all the nodes of the element.
After obtaining the new element values via Eqn.\eqref{eq:simple_avg} the point
averages are again evaluated via Eqn.\eqref{eq:volume_avg}. This process is repeated 
for a specified number of iterations (typically 1-5).
While very crude, this form of averaging works surprisingly well.

\subsection{$H^1$ (Weak) Laplacian Smoothing} \label{subsec:laplacian_smoothing}
In this case, the initial values $\alpha_0$ obtained for $\alpha$ 
are smoothed via:
\begin{equation} \label{eq:laplacian_smoothing1}
    \left[ 1 - \lambda \nabla^2 \right] \alpha = \alpha_0 ~~,~~
   \left. \alpha_{,n} \right|^{\Gamma} = 0 ~~.
\end{equation}

Here $\lambda$ is a free parameter which may be problem and mesh
dependent (its dimensional value is length squared). Discretization
via finite elements yields:
\begin{equation} \label{eq:laplacian_smoothing2}
    \left[ \Mmat_c + \lambda \Kmat_d \right] \alphavec = 
          \Mmat_{p1p0} \alphavec_0 ~~,
\end{equation}
where $\Mmat_c, \Kmat_d, \Mmat_{p1p0}$ denote the consistent mass 
matrix, the stiffness or `diffusion' matrix obtained for the Laplacian 
operator and the projection matrix from element values ($\alphavec_0$)
to point values ($\alphavec$).

\subsection{Pseudo-Laplacian Smoothing} \label{eq:pseudo_laplacian_smoothing}
One can avoid the dimensional dependency of $\lambda$ by smoothing via:
\begin{equation} \label{eq:pseudo_laplacian_smoothing1}
    \left[ 1 - \lambda \nabla h^2 \nabla \right] \alpha = \alpha_0 ~~,
\end{equation}
where $h$ is a characteristic element size. For linear elements, one
can show that this is equivalent to:
\begin{equation} \label{eq:pseudo_laplacian_smoothing2}
    \left[ \Mmat_c + \lambda \left(\Mmat_l - \Mmat_c \right) \right] 
      \alphavec = \Mmat_{p1p0} \alphavec_0 ~~,
\end{equation}
where $\Mmat_l$ denotes the lumped mass matrix \cite{lohner2008applied}. In the examples
shown below this form of smoothing was used for the gradients,
setting $\lambda=0.05$.

\section{Implementation in black-box solvers}
\label{sec:black_box_solvers}

The optimization cycle outlined above can be implemented in a very
efficient way if one has direct access to the source-code of
finite element-based structural mechanics solvers, but is also amenable 
to black-box (e.g. commercial) solvers. A possible way to proceed is
the following:
\begin{itemize}
\item[-] Output the original stiffness matrix $\Kmat_{el}$ for each element;
\item[-] For each optimization step/cycle:
\item[-] With the current element values for $\alpha$: build the
new stiffness matrix; this is usually done with a user-defined subroutine
or module (all commercial codes allow for that);
\item[-] For each load case $i$:
\item[-] Solve the forward problem ($\rightarrow \uvec_i$);
\item[-] Post-process the results of the forward problem in order
to obtain the displacements and strains;
\item[-] With the measured and computed displacements/
strains and weights: compute the cost function part of this load case
($\rightarrow I_i$);
\item[-] Build the `force vector' (i.e. the right-hand-side) for the
adjoint problem by comparing the measured and computed displacements/
strains and weighting them appropriately;
\item[-] Solve the adjoint problem ($\rightarrow \utildevec_i$);
\item[-] With $\uvec_i, \utildevec_i$ and the original stiffness matrix
$\Kmat_{el}$: get the gradient in each element;
\item[-] Smooth the gradients (either via a `fake-heat-solver' if
Laplacian smoothing is desired, or via an external smoother);
\item[-] Send the cost function and the smoothed gradients to the optimizer;
\item[-] Update $\alpha$
\end{itemize}

\section{Optimization of loadings}
\label{sec:optimization_loadings}

The aim of choosing a minimal yet optimal set of forces to monitor
the weakening of structures is to be able to determine the
field $\alpha(\xvec)$ as best as possible. From structural mechanics,
this implies that one should avoid regions where the strains are very
small or vanish. For these regions, $\alpha(\xvec)$ can assume an
arbitrary value without having any effect on the overall displacements
or strains. Therefore, {\bf the forces should be chosen such that the
number of regions with very small or vanishing strains should be
minimized}. \\
As stated before, for practical reasons the number and type of possible 
loads is limited. This `limitation of the search space for loads' opens 
up the possibility of a simple algorithm to determine the optimal 
choice. For each of the $n$ possible loads $\fvec_i, i=1,n$,
obtain the resulting displacement and strains, and record all
elements for which the strains are above a minimal sensor 
threshold $s_0$. \\
If only one force is to be applied, the obvious choice is to select
the one that produces the largest area with strains that are above
a minimal threshold. Having selected this force, the regions 
that have already been affected by this force (i.e. with strains 
that are above the minimal threshold) are excluded 
from further consideration. The next best force is then again the one that
is able to measure the largest area with strains that are above  
a minimal threshold. And so on recursively.

\section{Optimal placement of sensors}
\label{sec:optimal_placement}

Let us assume that a certain part $\Omega^w$ of the structure has 
weakened. This could be a region of several elements, or a single
element. This will lead to a change in the stiffness matrix, and
a resulting change in displacements and strains. The aim is to be
able to record and identify the spatial location of this weakening
with the minimum number of sensors. The change in displacements 
or strains due to a weakening requires the evaluation of the
derivative
\begin{equation} \label{eq:sensitivity1}
    D_{\alpha}(\xvec_i, \xvec_j) = {{\partial \uvec(\xvec_i) } 
                             \over {\partial \alpha(\xvec_j)}}
\end{equation}
for all possible combination of locations $\xvec_i, \xvec_j$.
In the most general case $\xvec_i$ is arbitrary, i.e. it could
be any node or element of the mesh. However, if sensors can only
be placed in certain regions of the domain, the location of
$\xvec_i$ can be reduced significantly. \\
Two possible ways were explored to obtain $ D_{\alpha}(\xvec_i, \xvec_j)$: 
forward-based and adjoint-based.

\subsection{Forward-based} \label{subsec:forward_based}
An immediate approach takes each `region of elements', changes
the stiffness matrix and computes the resulting changes in
displacements and strains at the possible sensor locations. Dropping the
index for the load cases, for each of them this results in:
\begin{equation} \label{eq:forward_based_1}
    ( \Kmat + \Delta \Kmat ) \cdot ( \uvec + \Delta \uvec ) = \fvec 
                                          ~~.
\end{equation}

With the original system (Eqn.\eqref{eq:forward}) this results in:
\begin{equation} \label{eq:forward_based_2}
    \Kmat \cdot \Delta \uvec = 
 - \Delta \Kmat \cdot ( \uvec + \Delta \uvec ) ~~.
\end{equation}

As the inverse (or LU decomposed) matrix of $\Kmat$ is assumed as
given (it was needed to compute $\uvec$), there are now two options: \\
\noi
a) Neglect the higher order terms $\Delta \Kmat \cdot \Delta \uvec$,
which results in:
\begin{equation} \label{eq:forward_based_3}
    \Kmat \cdot \Delta \uvec = - \Delta \Kmat \cdot \uvec 
                                                  ~~.
\end{equation}
\noi
b) Iterate for $\Delta \uvec$:
\begin{equation} \label{eq:forward_based_4}
    \Kmat \cdot \Delta \uvec^{i+1} = 
 - \Delta \Kmat \cdot ( \uvec + \Delta \uvec^{i} ) ~~,~~i=1,k ~~,~~
   \Delta \uvec^{0}=0
\end{equation}

One is now able to determine for each possible weakening region
$\Omega^w_k,~k=1,r$ the resulting deformations and strains, and
with the thresholds $u_0, s_0$ those sensors that are able to
monitor the weakening. \\
The effect of weakening a region (in the
extreme case a single element) on the sensors implies, in the
worst case, a CPU requirement that is of order
$O(N^2_{el} \cdot N_{bandwidth})$ for each of the $n$ load
cases, where $N_{el}$ denotes the
number of elements and $N_{bandwidth}$ the bandwidth of the system
matrix $\Kmat$. Clearly, for large problems with $N_{el}=O(10^6)$
this can become costly. Several options to manage this high cost
are treated in Section \ref{subsec:implementation_details}.

\subsection{Adjoint-based} \label{subsec:adjoint_based}
A more elegant (and faster) approach makes use of the adjoint
to obtain the desired sensitivities for each possible sensor
locations. The desired quantity whose derivative with respect
to element strength factor $\alpha$ is sought (e.g. displacement
at a sensor location) can be written as:
\begin{equation} \label{eq:adjoint_based_1}
    J = u(\text{bc},\text{loads},\alpha,\xvec) ~~
\end{equation}

The desired derivative is given by:
\begin{equation} \label{eq:adjoint_based_2}
    {{dJ } \over {d\alpha}} =
   {{\partial J } \over {\partial \alpha}} 
 + {{\partial J } \over {\partial u}} 
   {{\partial u } \over {\partial \alpha}}
                                                      ~~.
\end{equation}

This can be augmented to a Lagrangian by invoking the elasticity
equations, resulting in:
\begin{equation} \label{eq:adjoint_based_3}
    L^J = u(\text{bc},\text{loads},\alpha,\xvec) 
       + \utildevec \cdot (\Kmat \uvec - \fvec) ~~.
\end{equation}

The derivatives result in the usual systems of equations:
\begin{subequations}
\begin{eqnarray} \label{eq:adjoint_based_4}
    L^J_{,\utildevec} = \Kmat \uvec - \fvec = 0    ~~, \\
    L^J_{,\uvec} = {{\partial u(\text{bc},\text{loads},\alpha,\xvec)} \over
                   {\partial \uvec}} + \utildevec^T \Kmat = 0 
                                                        ~~, \\
    L^J_{,{\alpha_e}} = \utildevec^T \Kmat_e \uvec 
                                                        ~~.
\end{eqnarray}
\end{subequations}

Using the adjoint the information sought is evaluated in
the opposite order to the previous (forward, element-based)
procedure. While in the forward case an element/region was weakened
resulting in displacements/strains for all nodes/elements, in the 
adjoint case a location is selected and the effect of weakening 
each element on this location is obtained. \\
Observe that in this case the CPU requirement is of order
$O(m \cdot N_{el} \cdot N_{bandwidth})$ for each of the $n$
loadcases, where $m$ denotes the number of sensors (assumed to be 
much lower than the number of elements), $N_{el}$ the number of 
elements and $N_{bandwidth}$ the bandwidth of the system matrix $\Kmat$.
The advantages of using the adjoint may become even more pronounced
for nonlinear problems. While the forward-based procedure would
require the solution of a nonlinear problem for each element
(or element group), the adjoint always remains a linear problem.

\subsection{Sensor placement} \label{subsec:sensor_placement}
If only one sensor is to be placed, the obvious choice is to select
the one that is able to measure the highest number of weakening regions.
Having selected this sensor, the weakening regions that were able to
be measured are excluded from further consideration.
The next best sensor is then again the one that
is able to measure the highest number of the remaining weakening regions.
And so on recursively.

\subsection{Implementation details} \label{subsec:implementation_details}
There are two aspects of the sensitivity calculation procedures
for optimal sensor placement that need to be addressed: computation and
storage requirements.

\subsection{Computation} \label{subsec:computation}
The effect of weakening a region (in the extreme case a single element) 
for the forward option, or the sensitivity of each element/point in the
mesh for the adjoint option implies, in the worst case, a CPU requirement 
that is of order $O(N^2_{el} \cdot N_{bandwidth})$, where $N_{el}$ 
denotes the number of elements and $N_{bandwidth}$ the bandwidth of 
the system matrix $\Kmat$. Clearly, for large problems with 
$N_{el}=O(10^6)$ this can become costly. As is so often the case in
computational mechanics, algorithms and hardware can help alleviate
this problem.
\subsubsection{Clustering of elements}
Instead of weakening a single element (forward option) or computing 
the sensity of a single node/element (adjoint option), one can cluster
elements into subregions. The CPU requirements then
decrease to $O(N_{reg} \cdot N_{el} \cdot N_{bandwidth})$, where
$N_{reg}$ denotes the number of subregions. In the present case an
advancing front technique was used to cluster elements into subregions.
The size of the subregions can be specified via a minimum required number 
of elements per subregion, the area/volume of the subregion or the
minimum distance from the first element/point of the subregion.
Given that in 3-D the number of elements in a subregion increases 
quickly, the number of subregions can be substantially lower than the
number of elements, yielding a considerable reduction in CPU requirements.
\subsubsection{Parallel computing}
The matrix problem that needs to be solved to obtain the effect of 
weakening a region/element (forward option) or to obtain the sensitivity
of a region/element (adjoint option) is independent of other 
regions/elements, making the problem embarrassingly parallel.

\subsection{Storage} \label{subsec:storage}
Storing the effect of weakening a region (in the
extreme case a single element) or the sensitivity of all elements for
every element implies, in the worst case, a storage requirement that 
is of order $O(N^2_{el})$. Even if one is only interested in the
effect on  $m$ sensors this implies $O(m \cdot N_{el})$.
Clearly, for large problems with $N_{el}=O(10^6)$
and $m=O(10^4)$ this can become an issue. A simple way to diminish the
storage requirements is to store the on/off sensing in powers of 2:
\begin{equation} \label{eq:storage}
    s = \sum_{i=1}^m \kappa_i 2^i ~~
\end{equation}
where $\kappa_i$ is either $1$ or $0$ depending on whether the
sensor was activated or not.

\subsection{Sensor placement with regions} \label{subsec:sensor_placement_with_region}
In some cases, the weakening of all elements can be achieved with
only a few sensors (in the extreme case a single sensor). However,
placing a single sensor would preclude being able to precisely define 
weakening regions. Therefore, only the elements in the neighbourhood 
of the selected sensor are excluded from further consideration. As before,
the neighbourhood of the sensor can be specified via the number of
elements, the area/volume or the distance from the sensor. The remainder
of the procedure outlined above remains the same.


\section{Examples}
\label{sec:examples}

All the numerical examples were carried out using two finite
element codes. The first, FEELAST \cite{lohner2023feelast}, is a
finite element code based on simple linear (truss), triangular
(plate) and tetrahedral (volume) elements with constant material
properties per element that only solves the linear elasticity
equations. The second, CALCULIX \cite{dhondt2022calculix}, is a general, open source
finite element code for structural mechanical applications with
many element types, material models and options.
The optimization loops were steered via a simple shell-script for
the adjoint-based gradient descent method. In all cases, 
a `target' distribution of $\alpha(\xvec)$ was given, together with 
defined external forces $\fvec_{\Gamma}$.
The problem was then solved, i.e. the displacements $\uvec(\xvec)$
and strains $\svec(\xvec)$ were obtained and recorded at the
`measurement locations' $\xvec_j,~j=1,m$. This then yielded
the `measurement pair' $\fvec, \uvec_j,~j=1,m$ or
$\fvec, \svec_j,~j=1,m$ that was used to determine the material
strength distributions $\alpha(\xvec)$ in the field.

\subsection{Plate with hole} \label{subsec:plate_with_hole}
\begin{figure}[!hbt]
    \centering
    \includegraphics[width=0.78\textwidth]{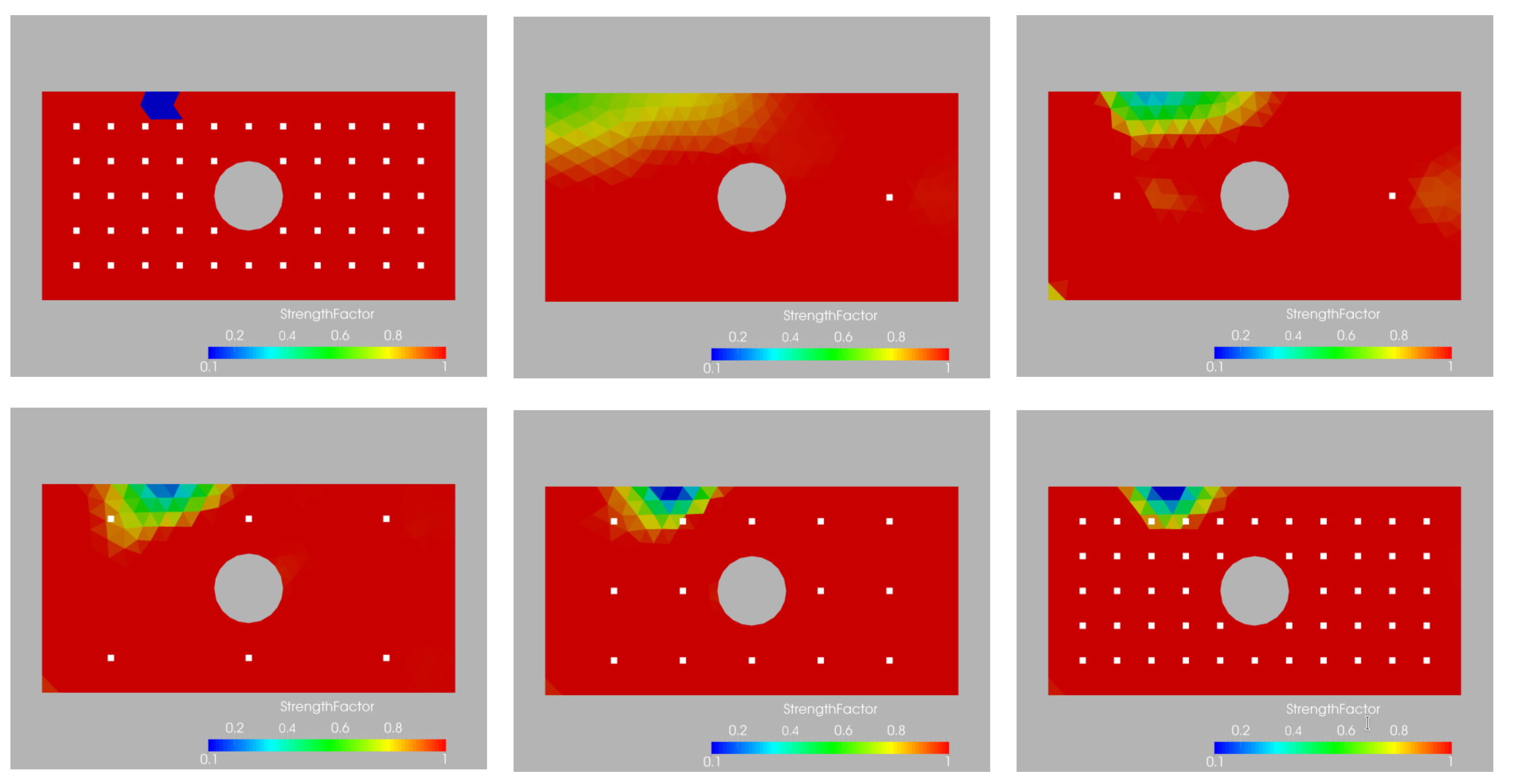}
    \caption{Plate With Hole: Effect of Sensors}
    \label{fig:plate_with_hole_effect}
\end{figure}

\begin{figure}[!hbt]
    \centering
    \includegraphics[width=0.78\textwidth]{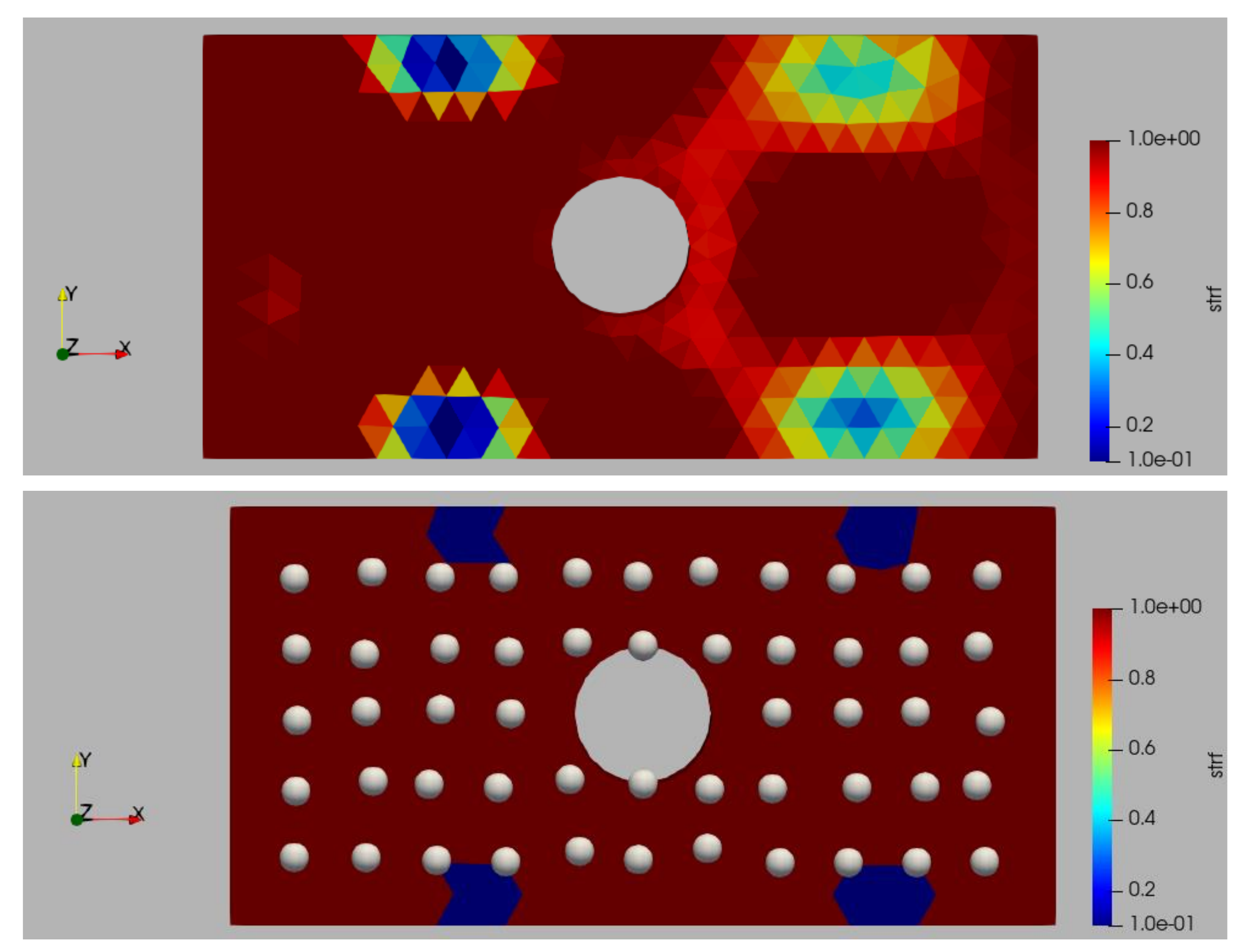}
    \caption{Plate With Hole: Sensing 4 Weakened Regions}
    \label{fig:plate_with_hole_4regions}
\end{figure}

The case is shown in Figures \ref{fig:plate_with_hole_effect} and \ref{fig:plate_with_hole_4regions}, and considers a plate 
with a hole. The plate dimensions are (all units in cgs):
$0 \le x \le 60$, $0 \le y \le 30$, $0 \le z \le 0.1$. A hole
of diameter $d=10$ is placed in the middle ($x=30, y=15$).
Density, Young's modulus and Poisson rate were set to
$\rho=7.8, E=2 \cdot 10^{12}, \nu=0.3$ respectively.
672 linear, triangular, plain stress elements were used.
The left boundary of the plate is assumed clamped ($\uvec=0$),
while a horizontal load of $q=(10^5,0,0)$ was prescribed at the
right end. In the first instance, a small region of weakened
material was specified. The ability of the procedure to detect
or `recover' this weakening based on the number of sensors used
(shown as white dots) is clearly visible. As the number of 
sensors increases, the region is recovered. Notice that even 
with~1 sensor a weakening is already detected, and that 
with~6 sensors the weakened region is clearly defined.

\ms \noi
In the second case, four weakening regions were specified.

\subsection{Thick plate with conical hole} \label{subsec:thick_plate}
\begin{figure}[!hbt]
    \centering
    \includegraphics[width=0.78\textwidth]{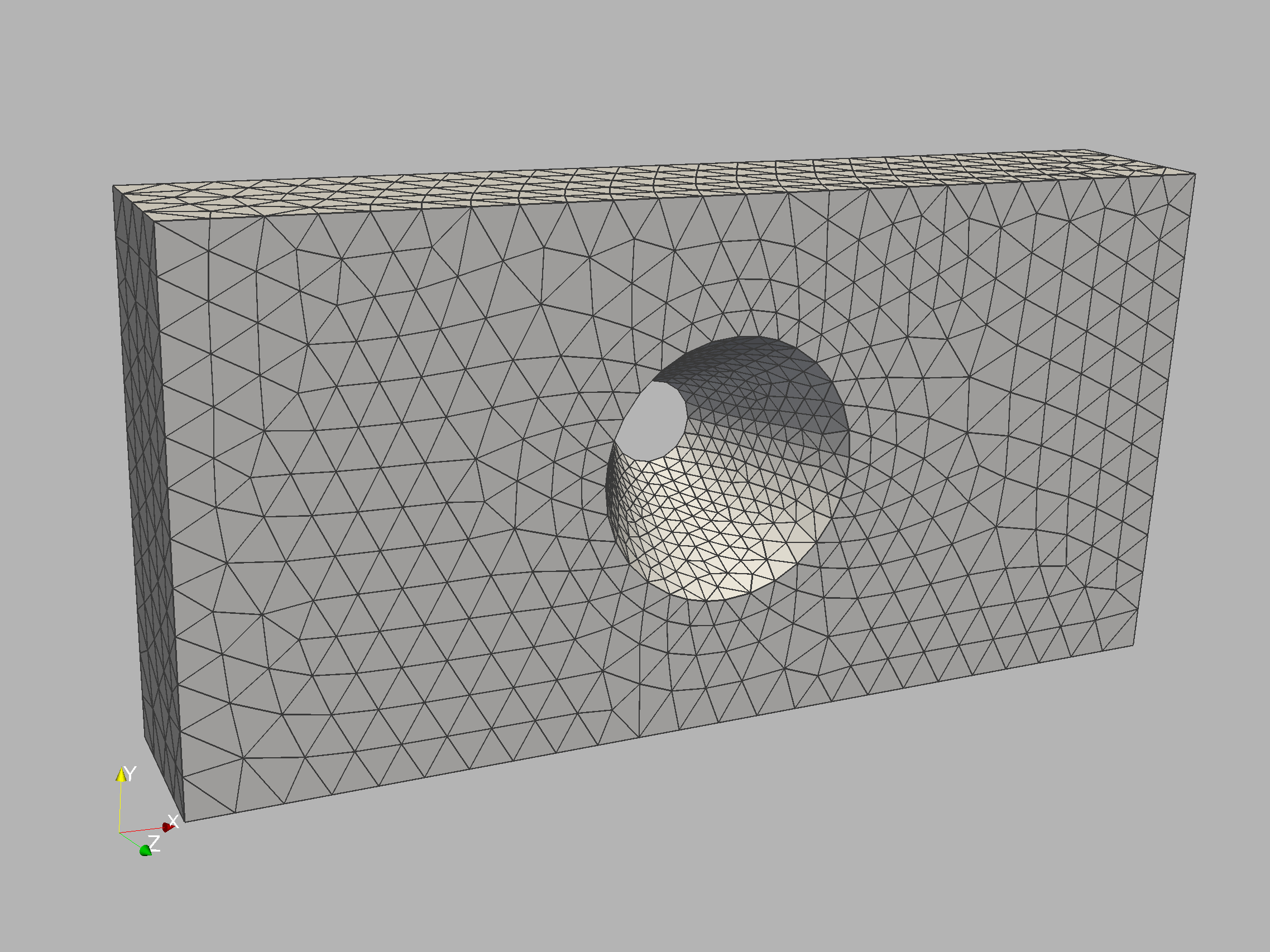}
    \caption{Thick Plate With Conical Hole: Surface of Coarse Mesh}
    \label{fig:thick_plate_coarse_mesh}
\end{figure}

\begin{figure}[!hbt]
    \centering
    \includegraphics[width=0.78\textwidth]{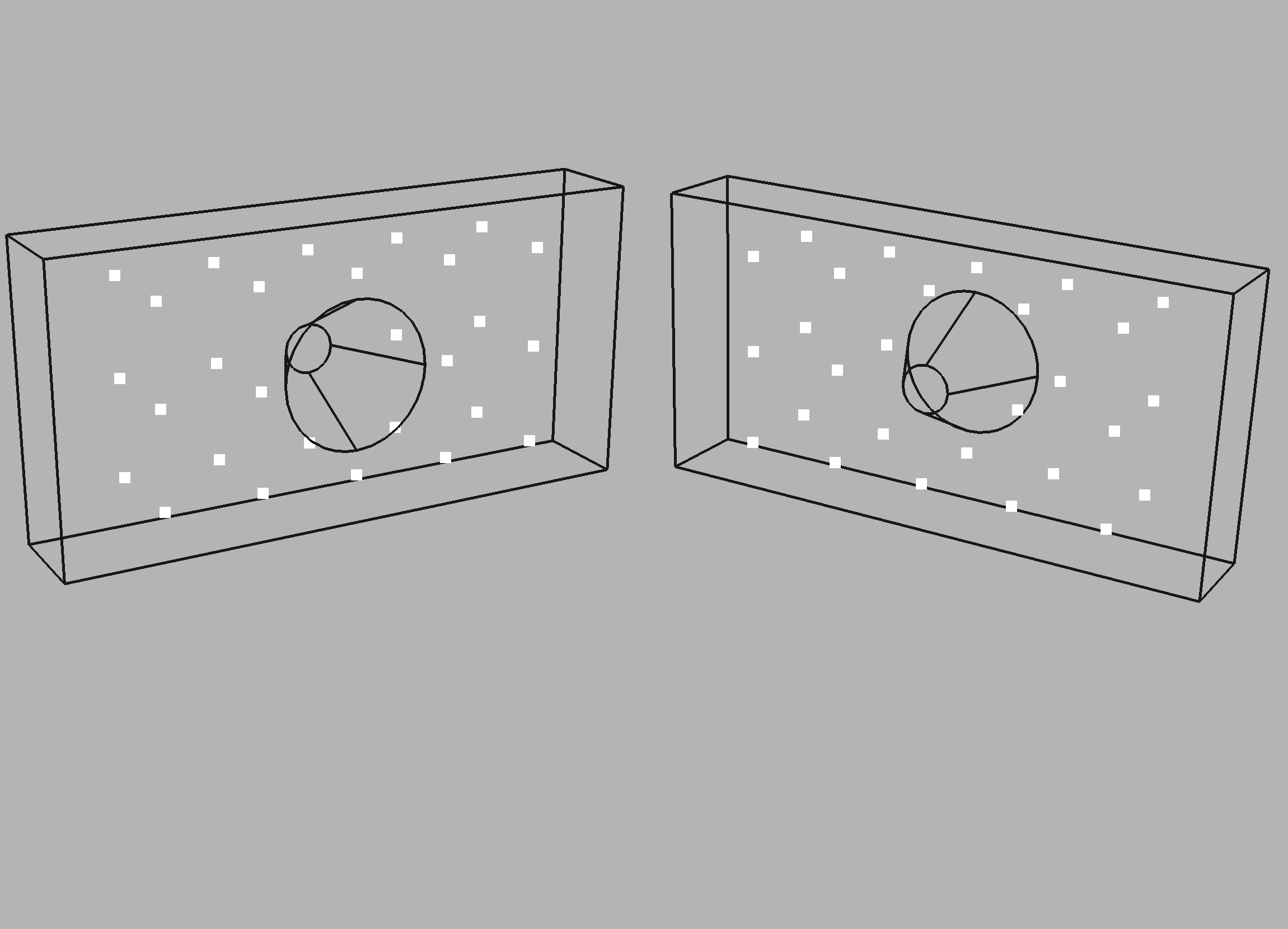}
    \caption{Thick Plate With Conical Hole: Location of Sensors}
    \label{fig:thick_plate_sensor_location}
\end{figure}

\begin{figure}[!hbt]
    \centering
    \includegraphics[width=0.78\textwidth]{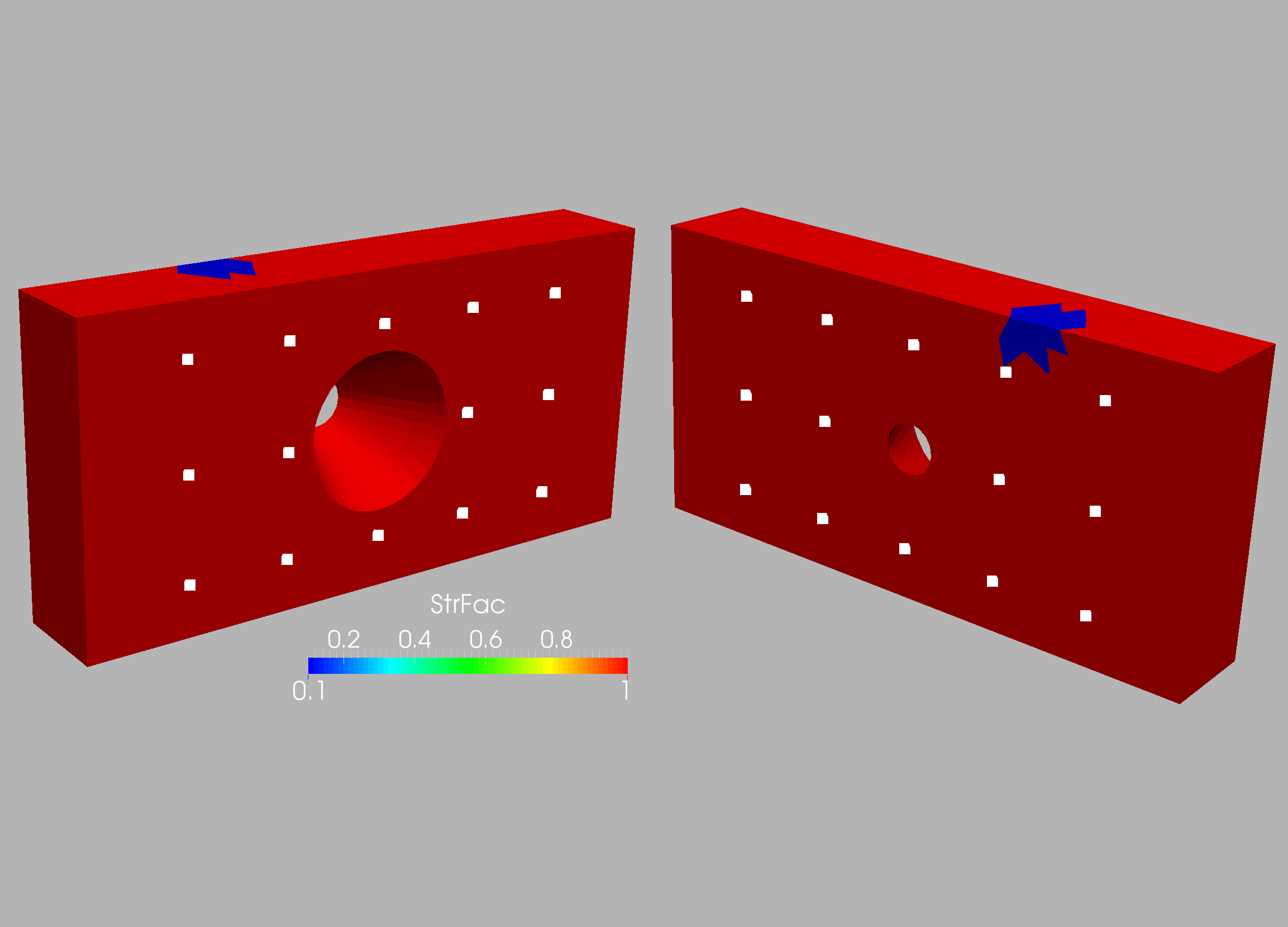}
    \caption{Thick Plate With Conical Hole: Weakened Region (Coarse Mesh)}
    \label{fig:thick_plate_target_coarse}
\end{figure}

\begin{figure}[!hbt]
    \centering
    \includegraphics[width=0.78\textwidth]{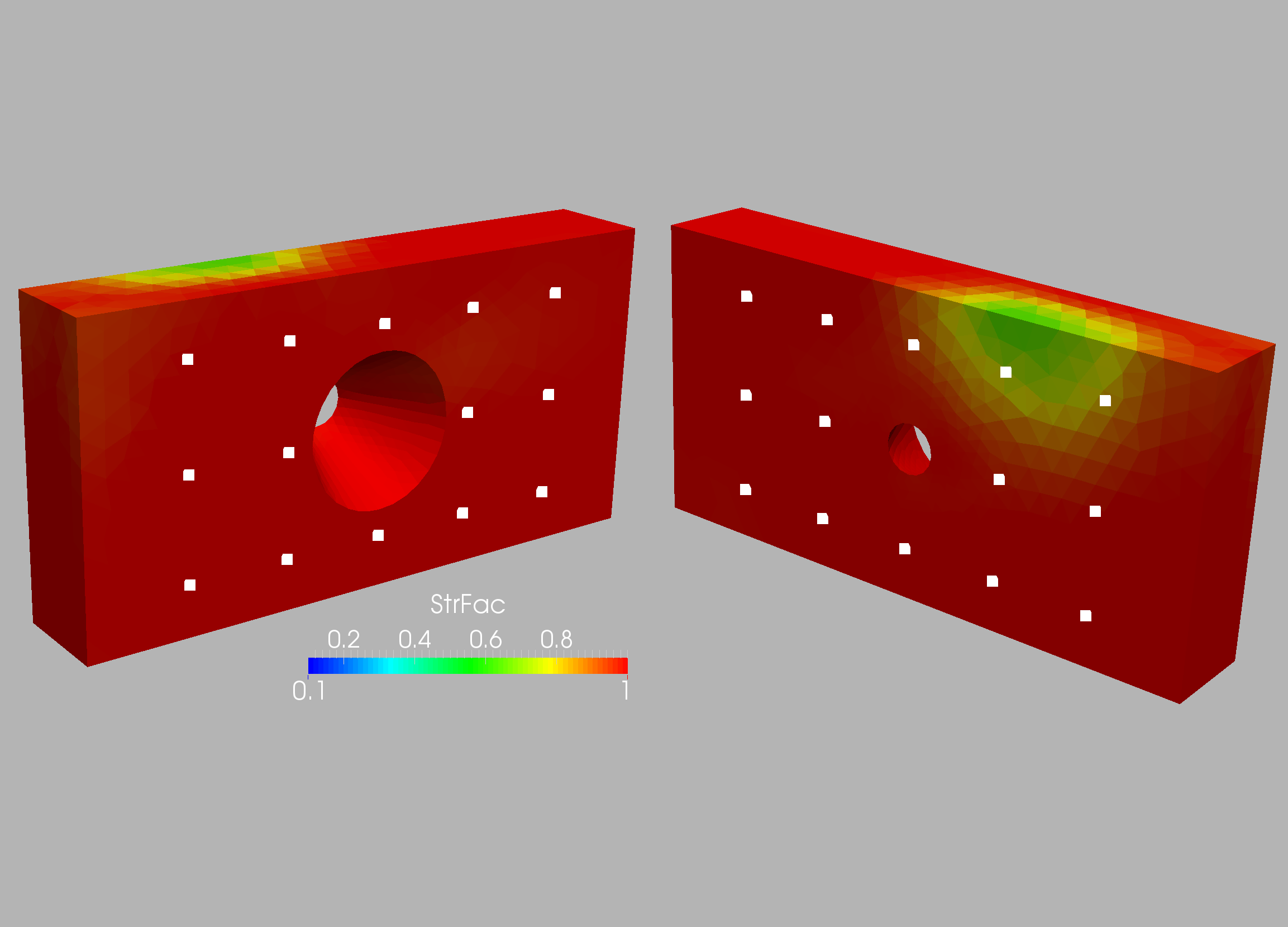}
    \caption{Thick Plate With Conical Hole: Computed Strength Factor With 28 Sensors (Coarse Mesh)}
    \label{fig:thick_plate_strf_coarse}
\end{figure}

\begin{figure}[!hbt]
    \centering
    \includegraphics[width=0.78\textwidth]{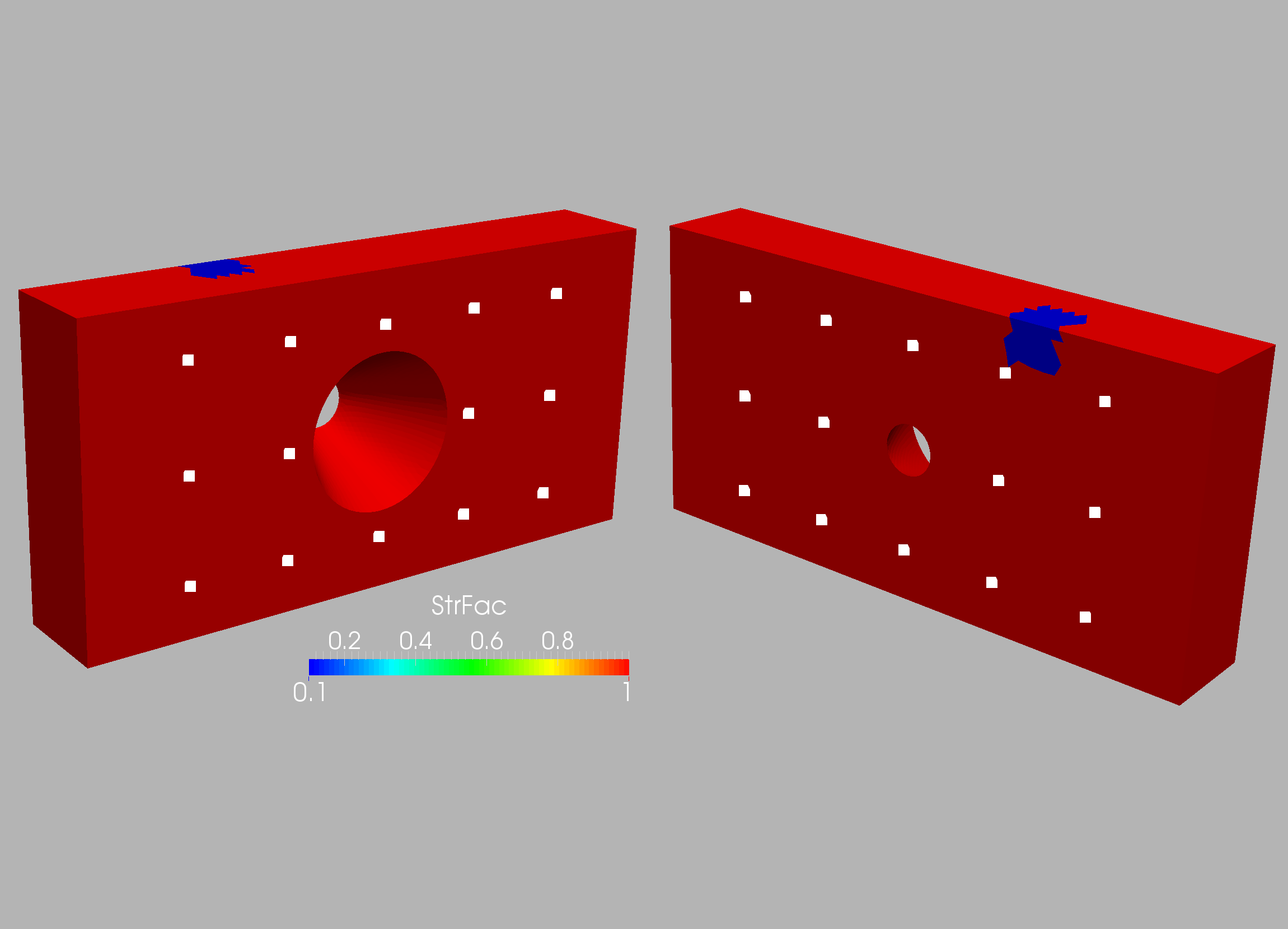}
    \caption{Thick Plate With Conical Hole: Weakened Region (Fine Mesh)}
    \label{fig:thick_plate_target_fine}
\end{figure}

\begin{figure}[!hbt]
    \centering
    \includegraphics[width=0.78\textwidth]{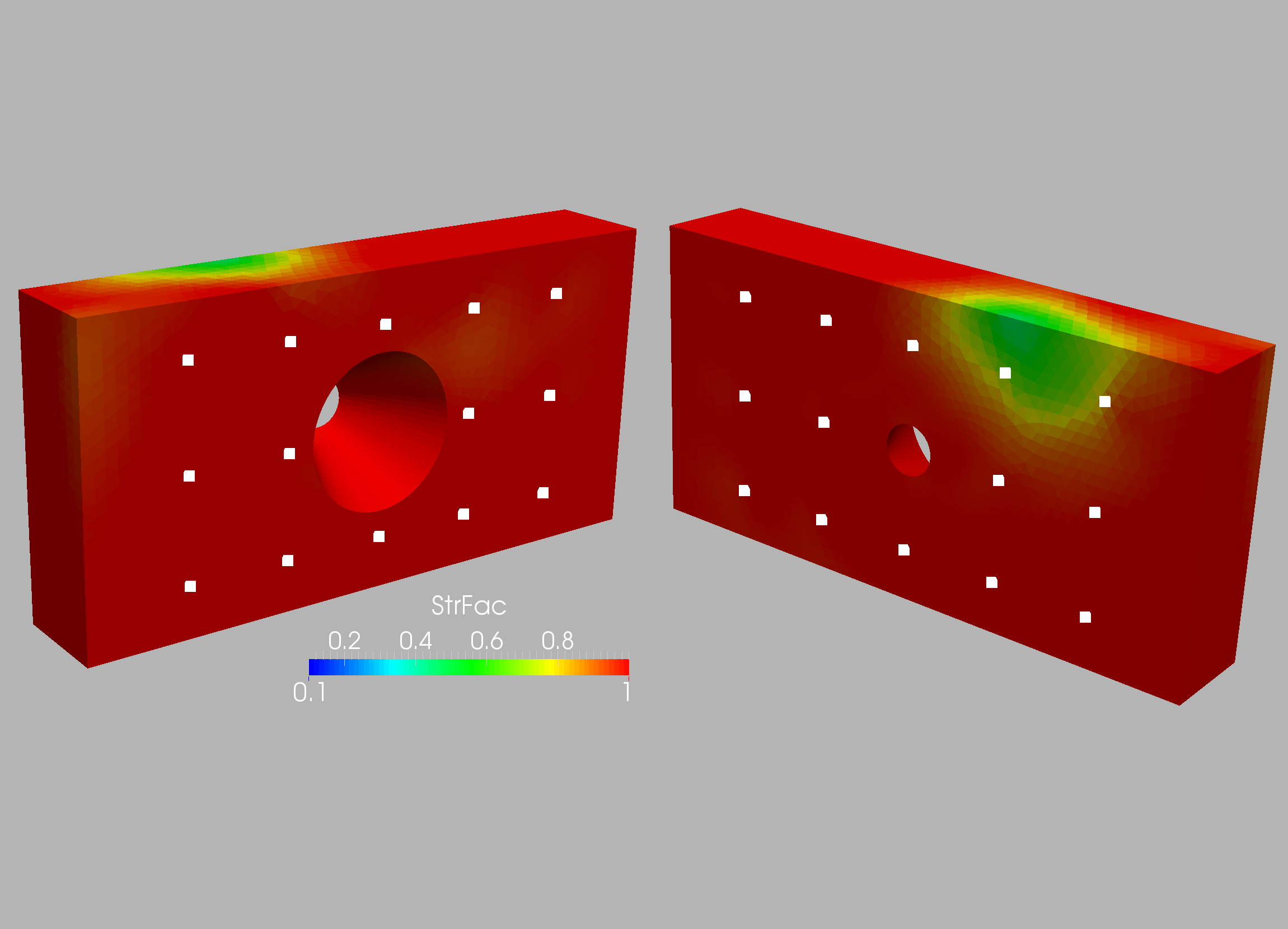}
    \caption{Thick Plate With Conical Hole: Computed Strength Factor With 28 Sensors (Fine Mesh)}
    \label{fig:thick_plate_strf_fine}
\end{figure}

\begin{figure}[!hbt]
    \centering
    \includegraphics[width=0.78\textwidth]{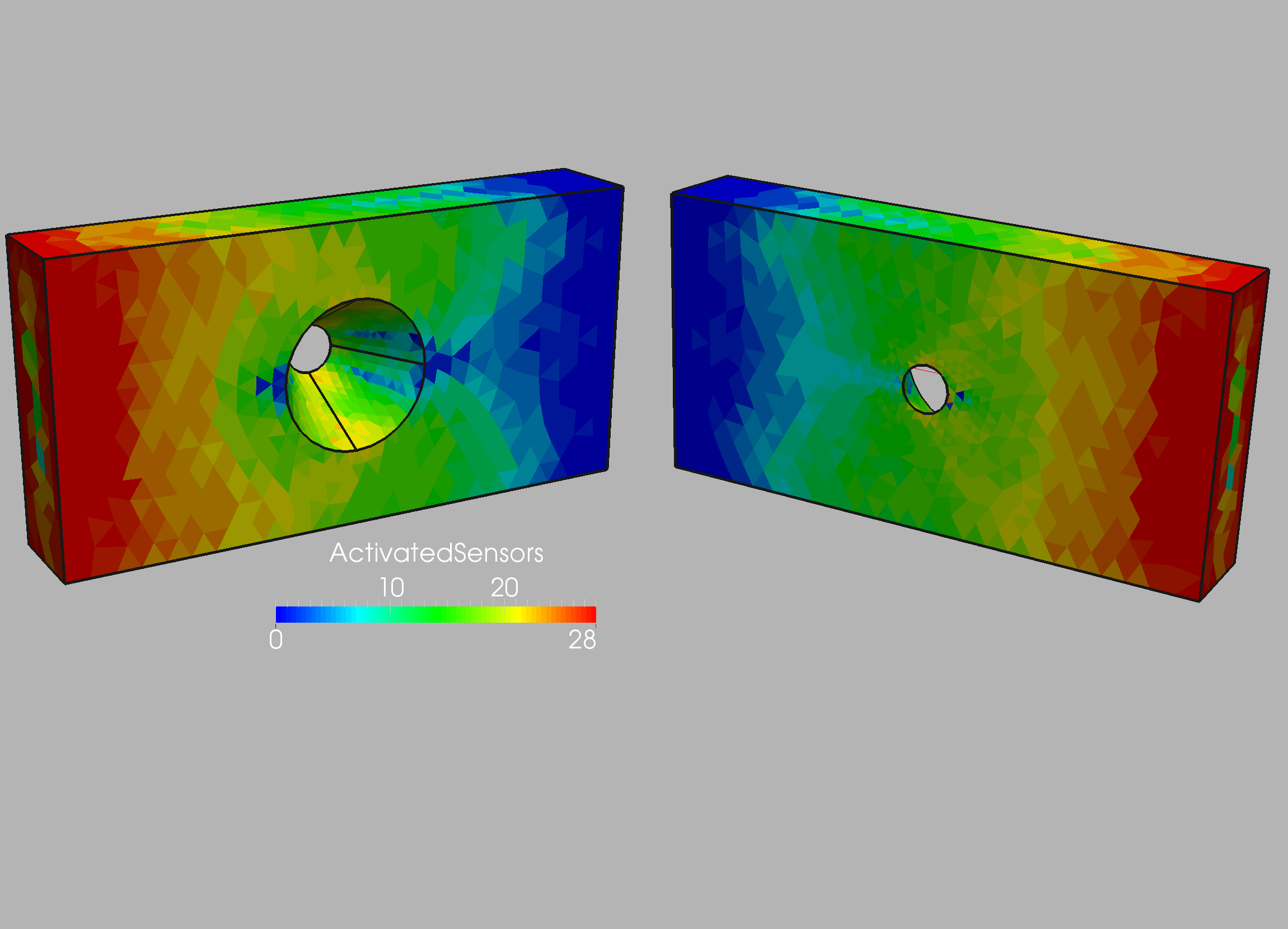}
    \caption{Thick Plate With Conical Hole: Nr. of Sensors Activated by Weakening An Element}
    \label{fig:thick_plate_activated_sensors}
\end{figure}

\begin{figure}[!hbt]
    \centering
    \includegraphics[width=0.78\textwidth]{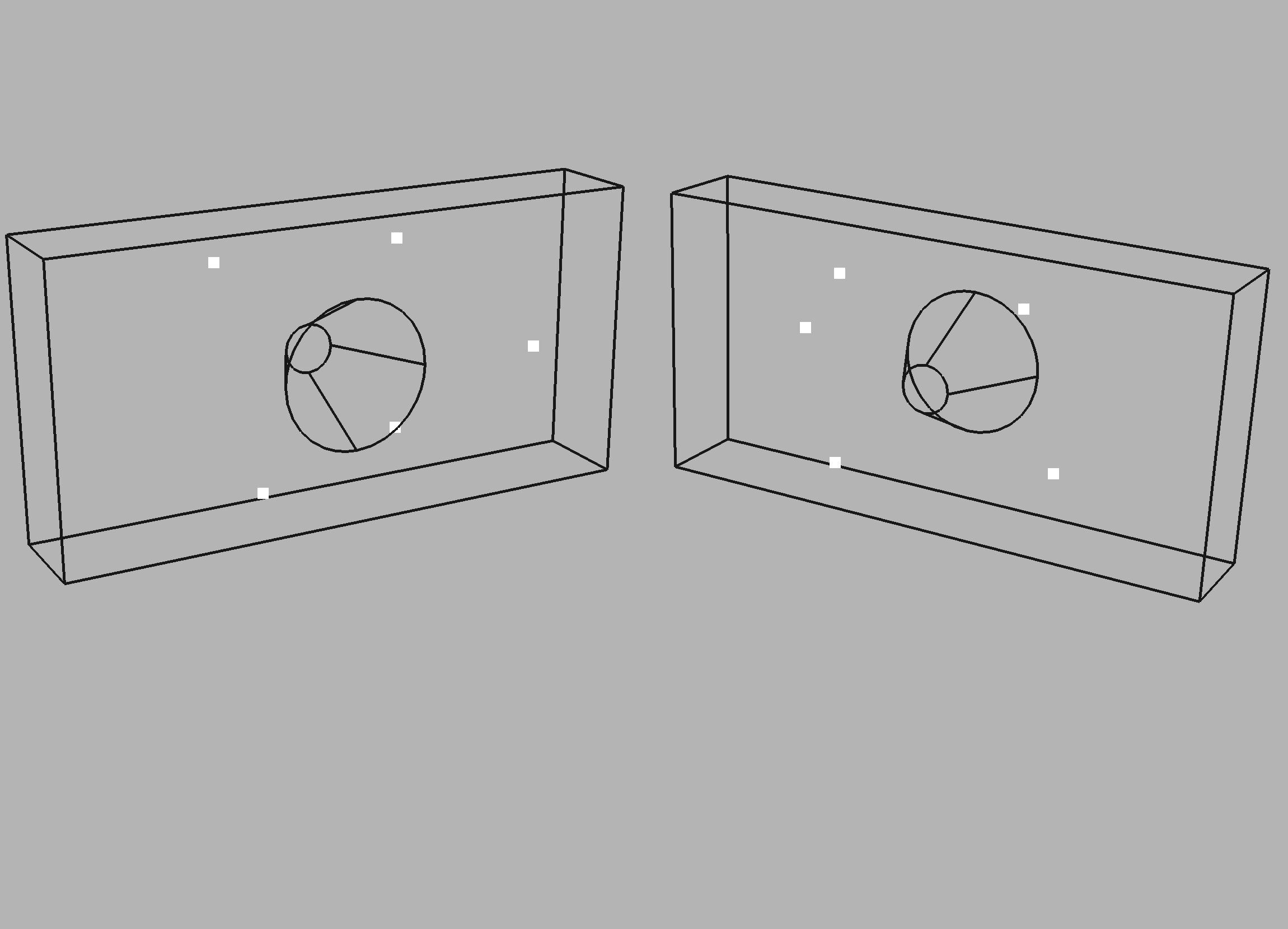}
    \caption{Thick Plate With Conical Hole: Optimal Location of Displacement Sensors}
    \label{fig:thick_plate_best_pts}
\end{figure}

\begin{figure}[!hbt]
    \centering
    \includegraphics[width=0.78\textwidth]{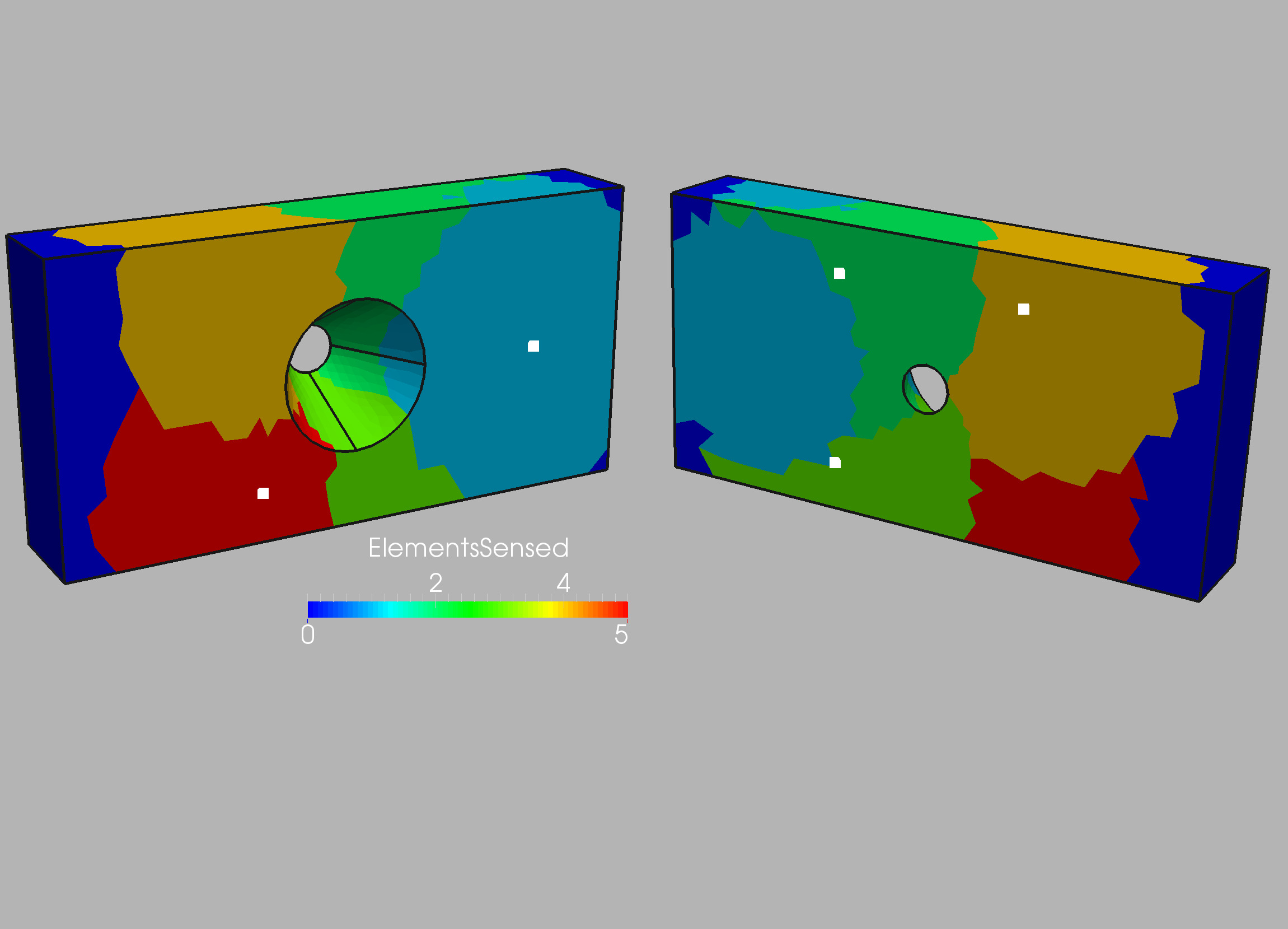}
    \caption{Thick Plate With Conical Hole: Elements `Sensed' for Each Displacement Sensor}
    \label{fig:thick_plate_elements_sensed}
\end{figure}

\begin{figure}[!hbt]
    \centering
    \includegraphics[width=0.78\textwidth]{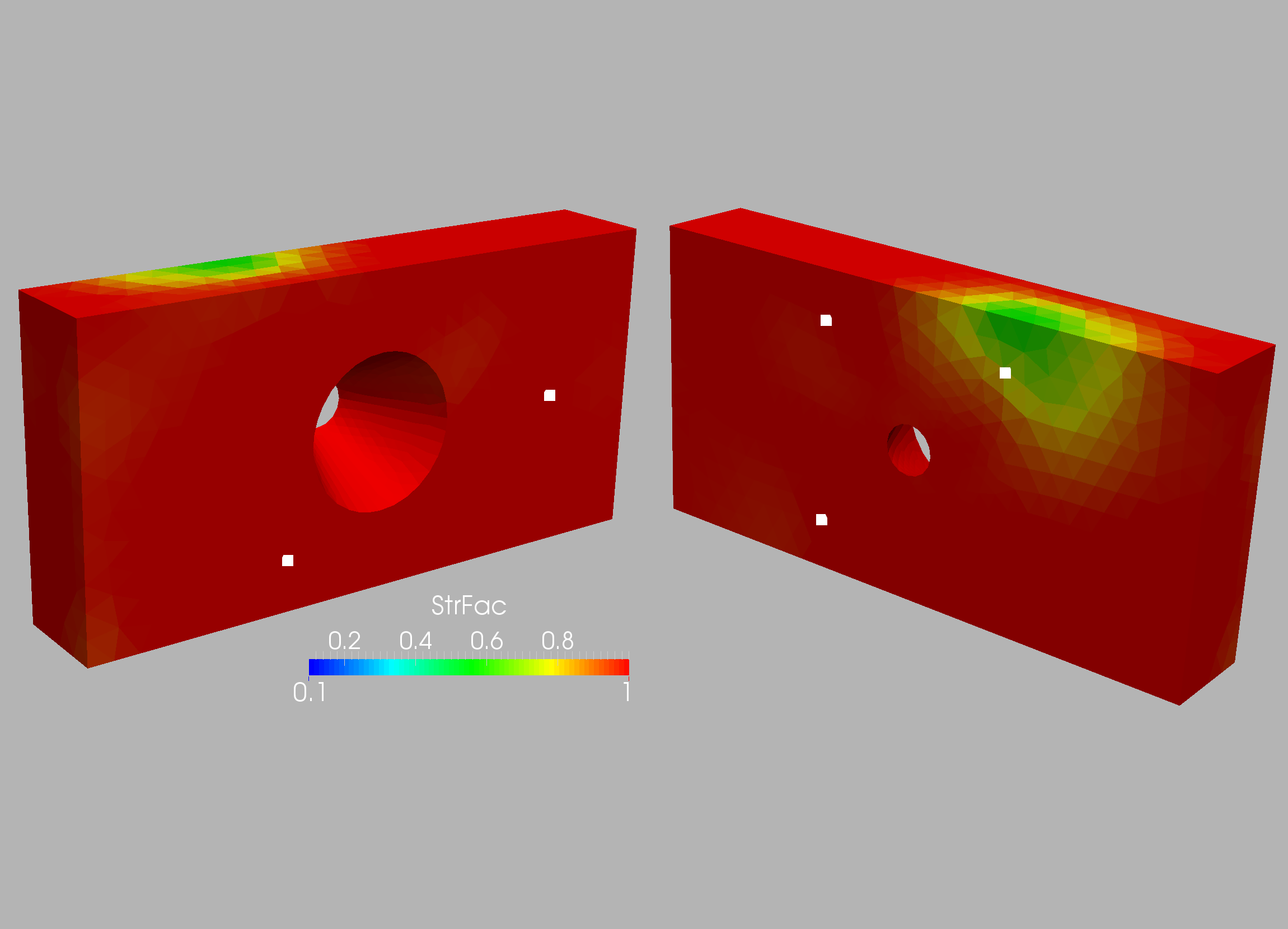}
    \caption{Thick Plate With Conical Hole: Computed Strength Factor With 5 Sensors (Coarse Mesh)}
    \label{fig:thick_plate_strf_coarse_5sensor}
\end{figure}

\begin{figure}[!hbt]
    \centering
    \includegraphics[width=0.78\textwidth]{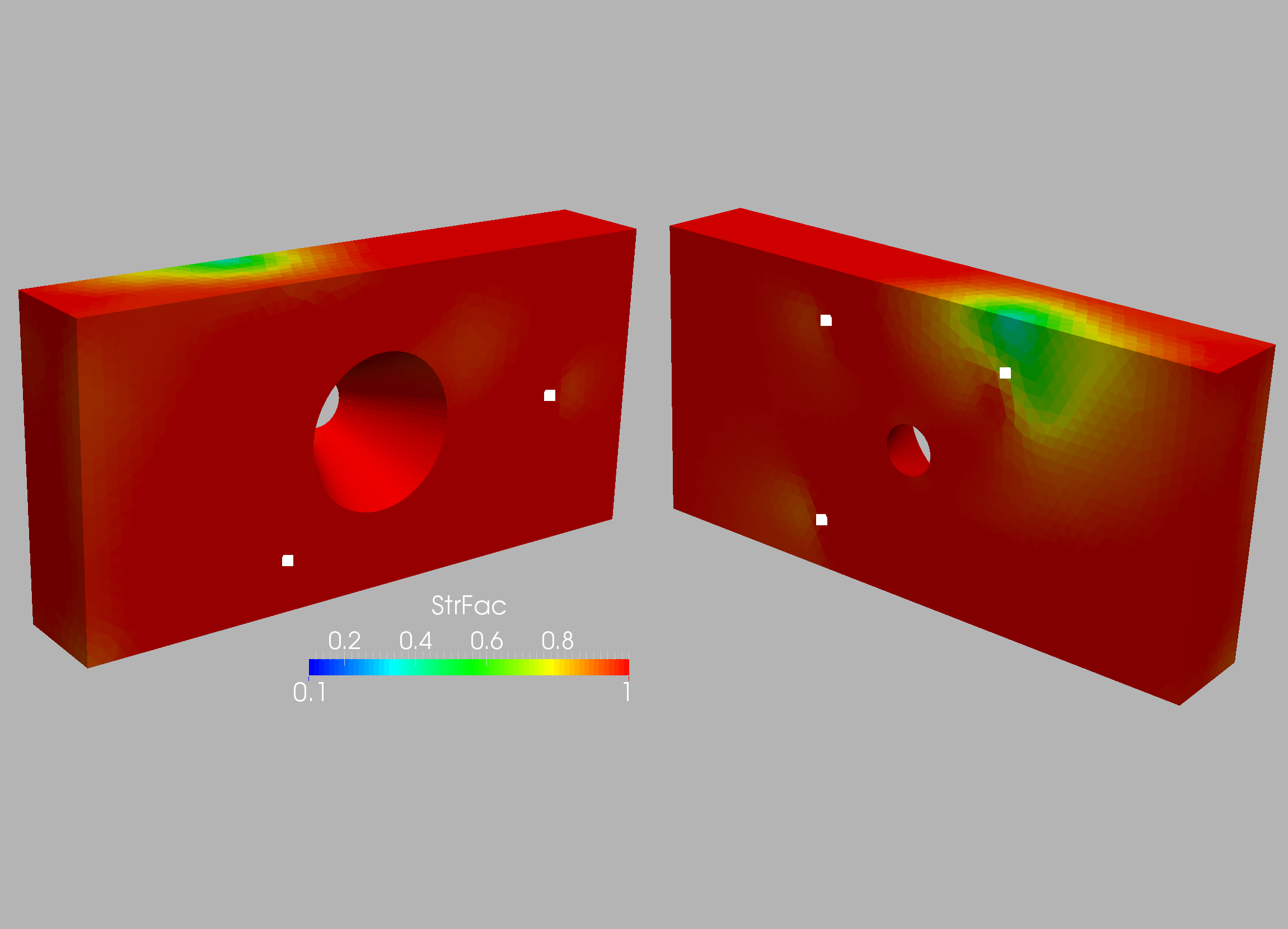}
    \caption{Thick Plate With Conical Hole: Computed Strength Factor With 5 Sensors (Fine Mesh)}
    \label{fig:thick_plate_strf_fine_5sensor}
\end{figure}

\begin{figure}[!hbt]
    \centering
    \includegraphics[width=0.78\textwidth]{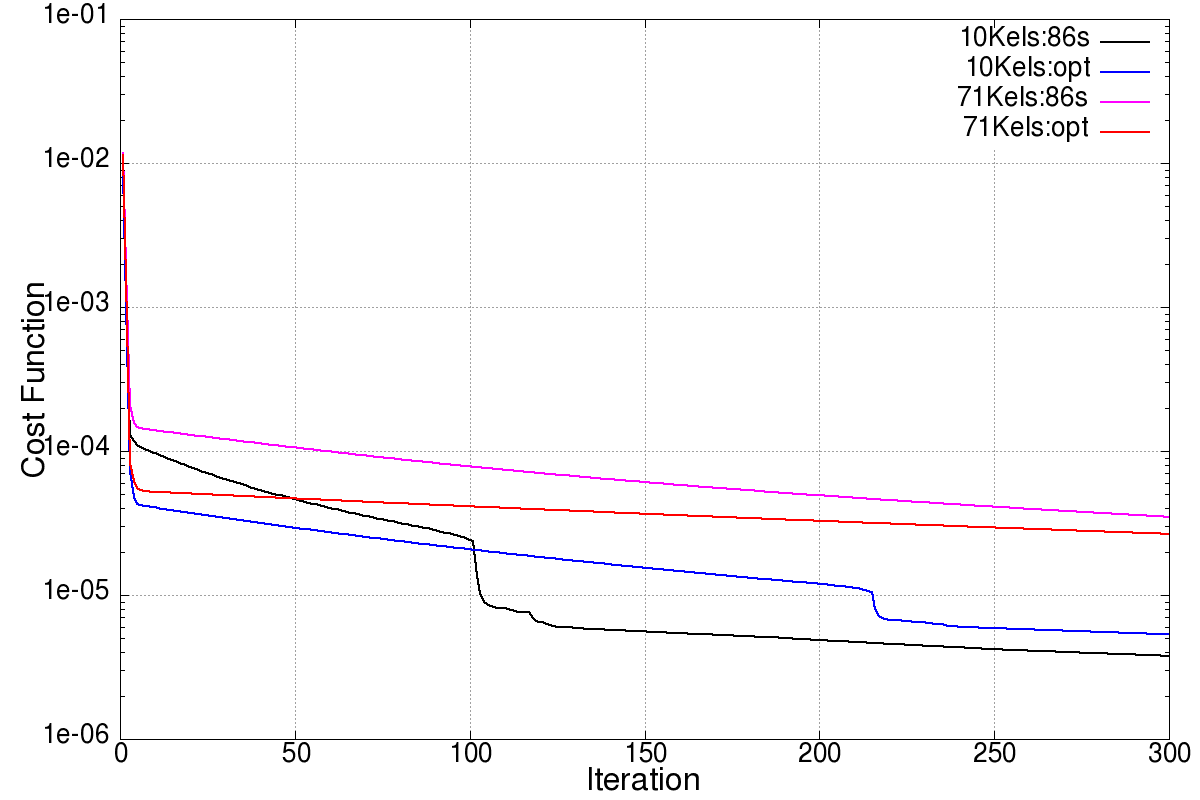}
    \caption{Thick Plate With Conical Hole: Convergence of Cost Function for all Cases}
    \label{fig:thick_plate_convergence}
\end{figure}

\begin{figure}[!hbt]
    \centering
    \includegraphics[width=0.78\textwidth]{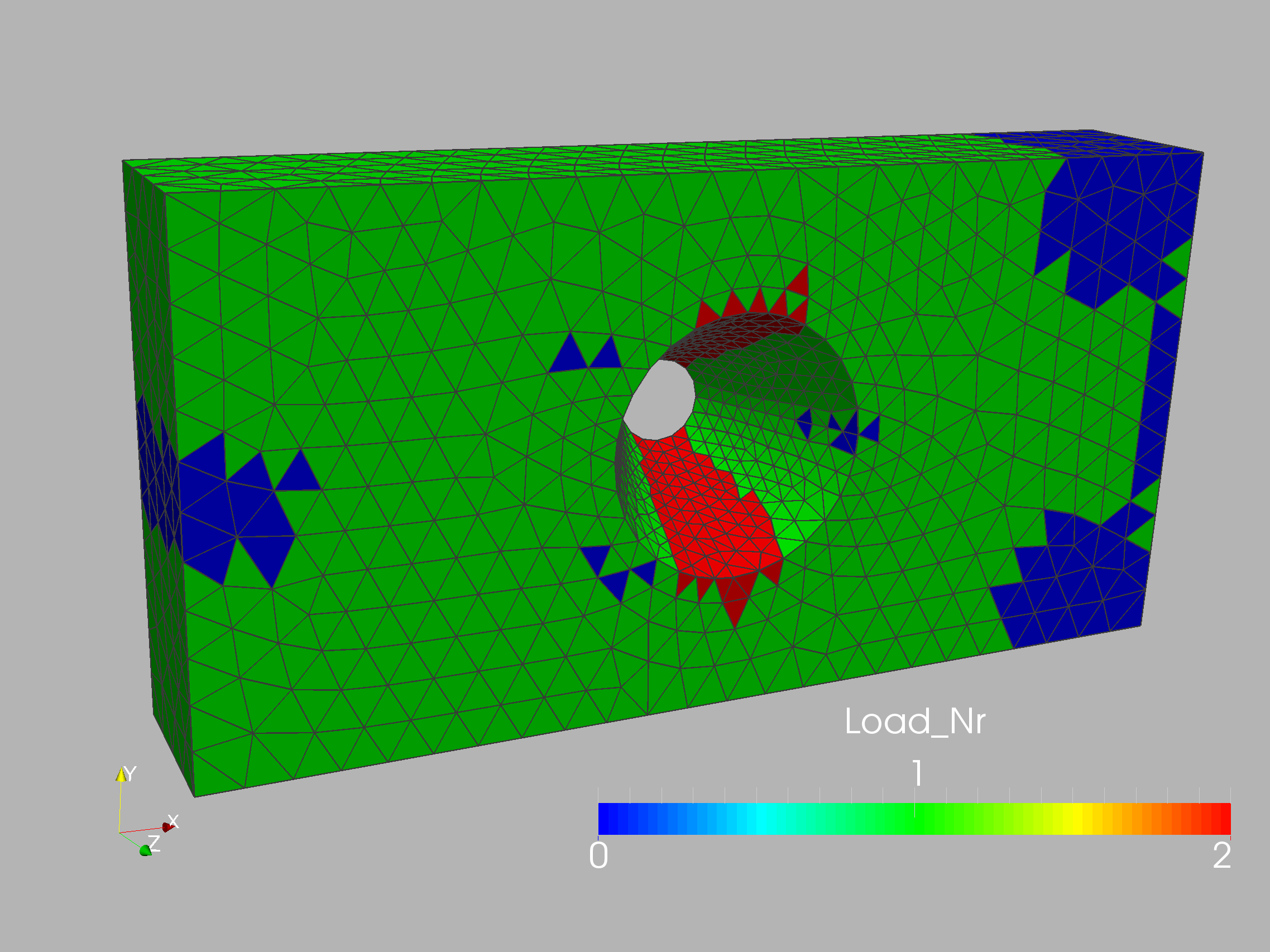}
    \caption{Thick Plate With Conical Hole: Loads Affecting Strain in Elements}
    \label{fig:thick_plate_optiload}
\end{figure}

The case is shown in Figure \ref{fig:thick_plate_coarse_mesh} and considers a thick plate
with a conical hole. The plate dimensions are (all units in cgs):
$0 \le x \le 60$, $0 \le y \le 30$, $0 \le z \le 10$. A conical hole
of diameter $d_1=5$ and $d_2=15$ is placed in the middle ($x=30, y=15$).
Density, Young's modulus and Poisson rate were set to
$\rho=7.8, E=2 \cdot 10^{12}, \nu=0.3$ respectively.
Two grids, of 19~K and 120~K linear tetrahedral elements (tet) were
used. The surface mesh of the coarser mesh
is shown in Figure \ref{fig:thick_plate_coarse_mesh}. A first series of runs with the 28~sensors
shown in Figure \ref{fig:thick_plate_sensor_location} were conducted. The target and computed weakening
for the two grids and the 28~sensors are shown in Figures \ref{fig:thick_plate_target_coarse}-\ref{fig:thick_plate_strf_fine}.
Note the proper detection of the weakened region. \\
Having proven that the technique works, the optimal number of
sensors and loads were obtained. Figure \ref{fig:thick_plate_activated_sensors} records the number
of (displacement) sensors that were able to measure/`sense' the
weakening of each element using the `forward-based' approach. 
As expected, the number is higher for
the elements close to the clamped boundary. The technique outlined
above (Section \ref{subsec:sensor_placement}) sorted the sensors. This resulted in just 5~sensors, 
whose location and `zone of influence' is shown in Figures~\ref{fig:thick_plate_best_pts} and \ref{fig:thick_plate_elements_sensed}. 
The weakening regions computed with these sensors for the two grids
are shown in Figures \ref{fig:thick_plate_strf_coarse_5sensor} and \ref{fig:thick_plate_strf_fine_5sensor}. One can see that even with this small number
of sensors the regions are well defined. The convergence history of the
cost function for these cases is plotted in Figure \ref{fig:thick_plate_convergence}.
Finally, Figure \ref{fig:thick_plate_optiload} shows the number of elements with strains above
a threshold for the 2 load cases: $q_1=(10^5,0,0)$ and $q_2=(0,-10^5,0)$.
One can see that with these 2 loads cases almost all elements are
affected, so any further load cases would not lead to more information.

%


\subsection{Connecting rod} \label{subsec:connecting_rod}

\begin{figure}[!hbt]
    \centering
    \includegraphics[width=0.78\textwidth]{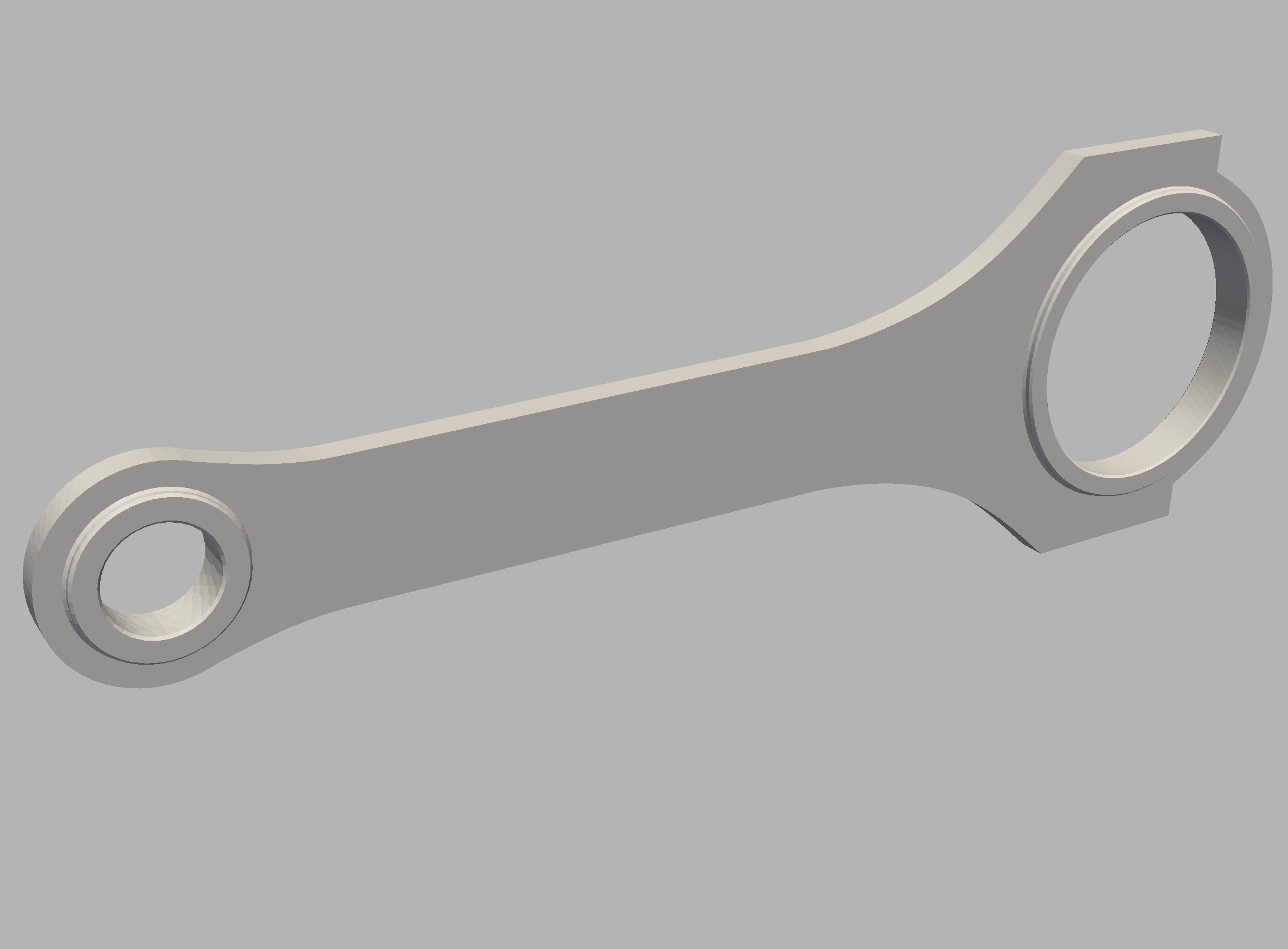}
    \caption{Connecting Rod: Surface}
    \label{fig:connecting_rod_surface}
\end{figure}

\begin{figure}[!hbt]
    \centering
    \includegraphics[width=0.78\textwidth]{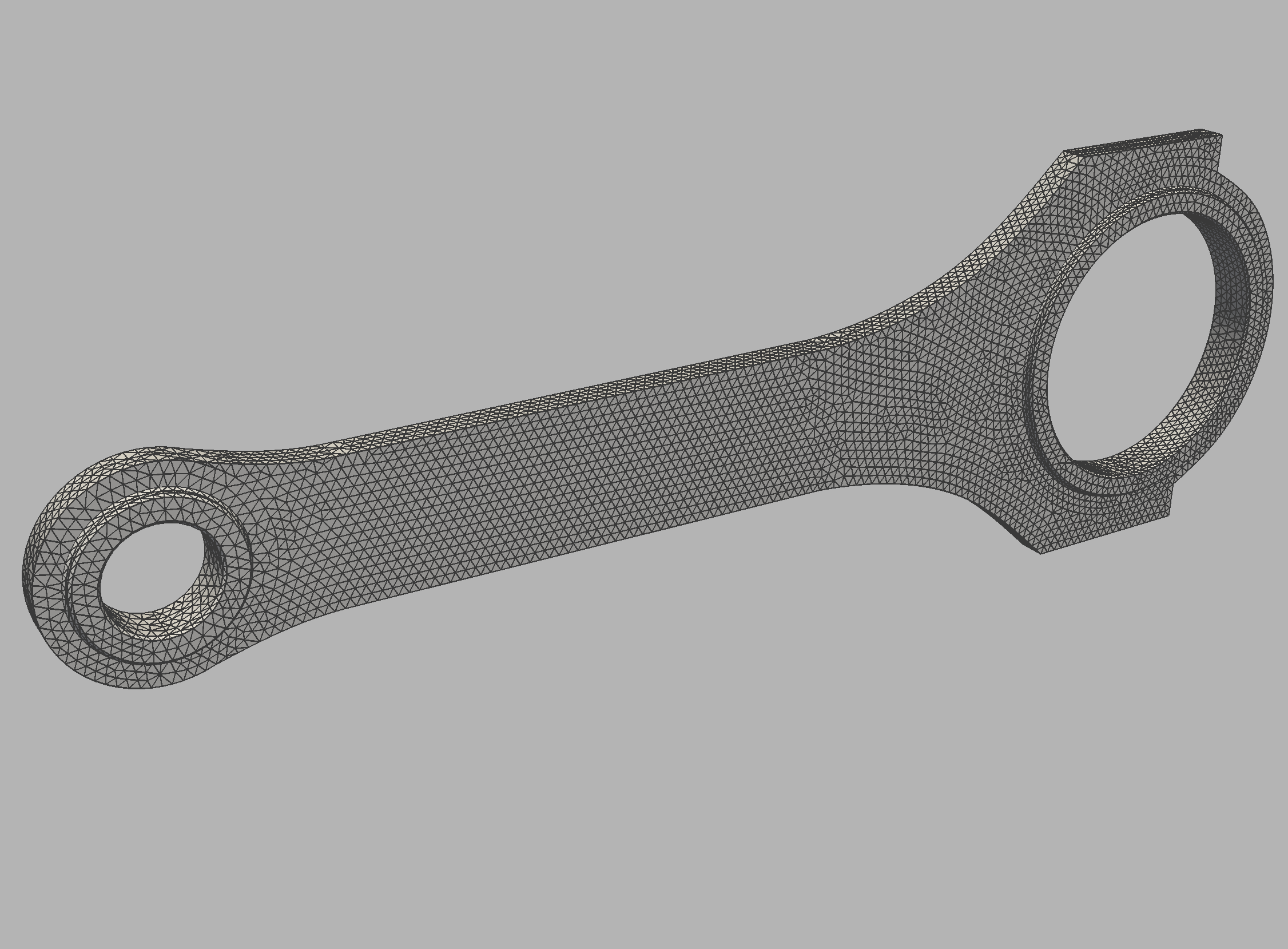}
    \caption{Connecting Rod: Surface Triangulation}
    \label{fig:connecting_rod_wire}
\end{figure}

\begin{figure}[!hbt]
    \centering
    \includegraphics[width=0.78\textwidth]{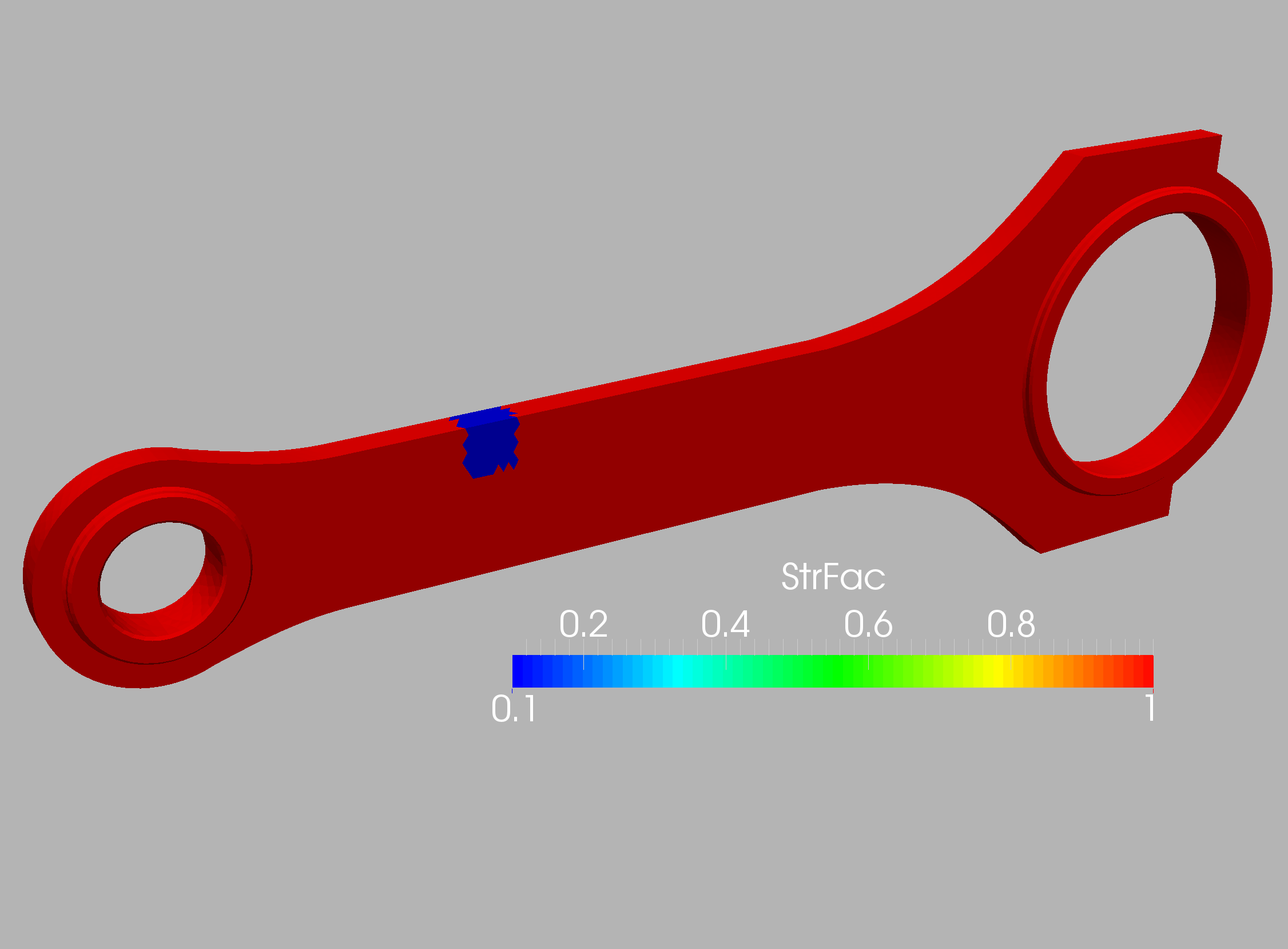}
    \caption{Connecting Rod: Target Strength Factor}
    \label{fig:connecting_rod_target_strf}
\end{figure}

\begin{figure}[!hbt]
    \centering
    \includegraphics[width=0.78\textwidth]{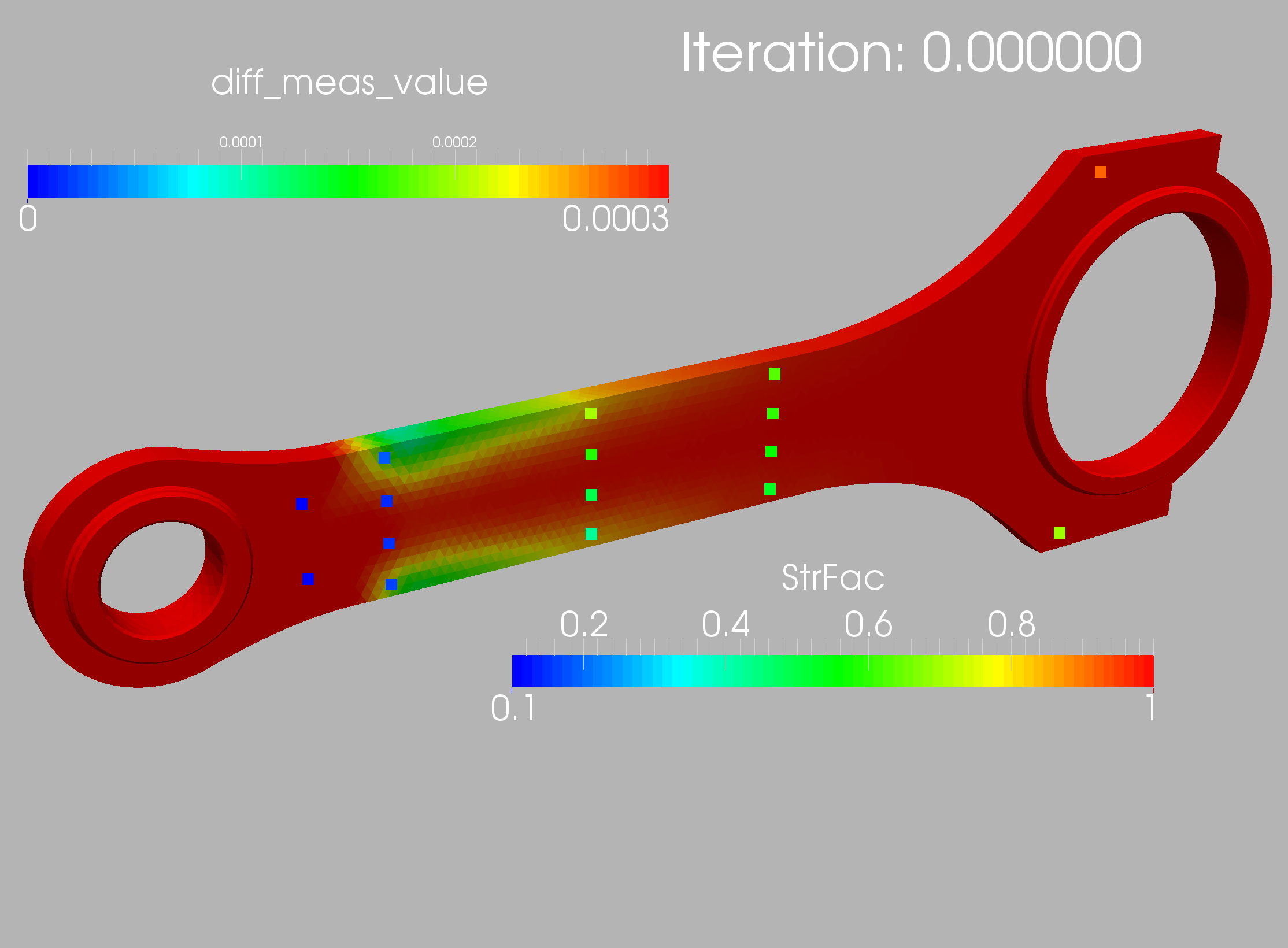}
    \caption{Connecting Rod: Iteration 0}
    \label{fig:connecting_rod_iter0}
\end{figure}

\begin{figure}[!hbt]
    \centering
    \includegraphics[width=0.78\textwidth]{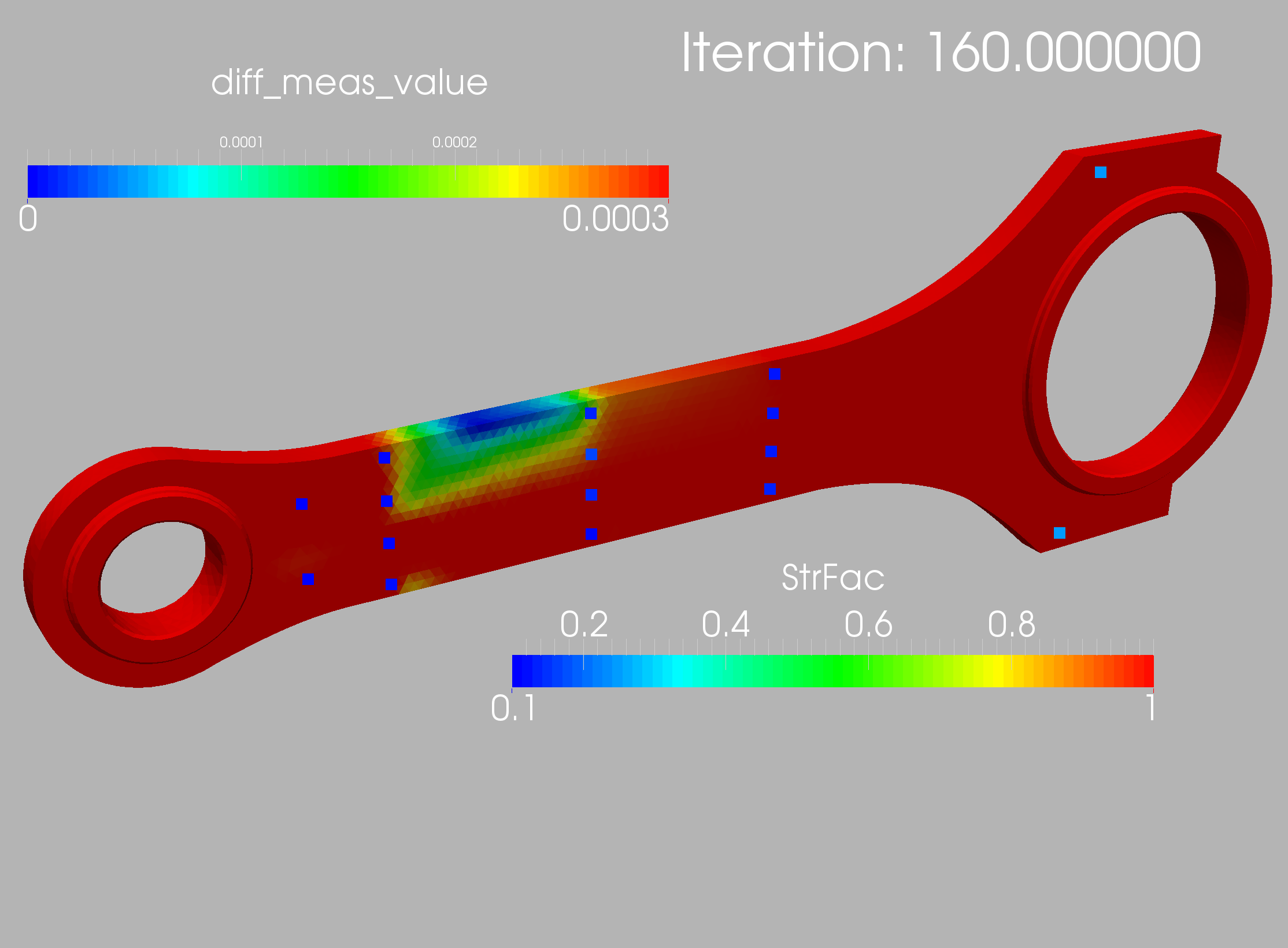}
    \caption{Connecting Rod: Iteration 160}
    \label{fig:connecting_rod_iter160}
\end{figure}

\begin{figure}[!hbt]
    \centering
    \includegraphics[width=0.78\textwidth]{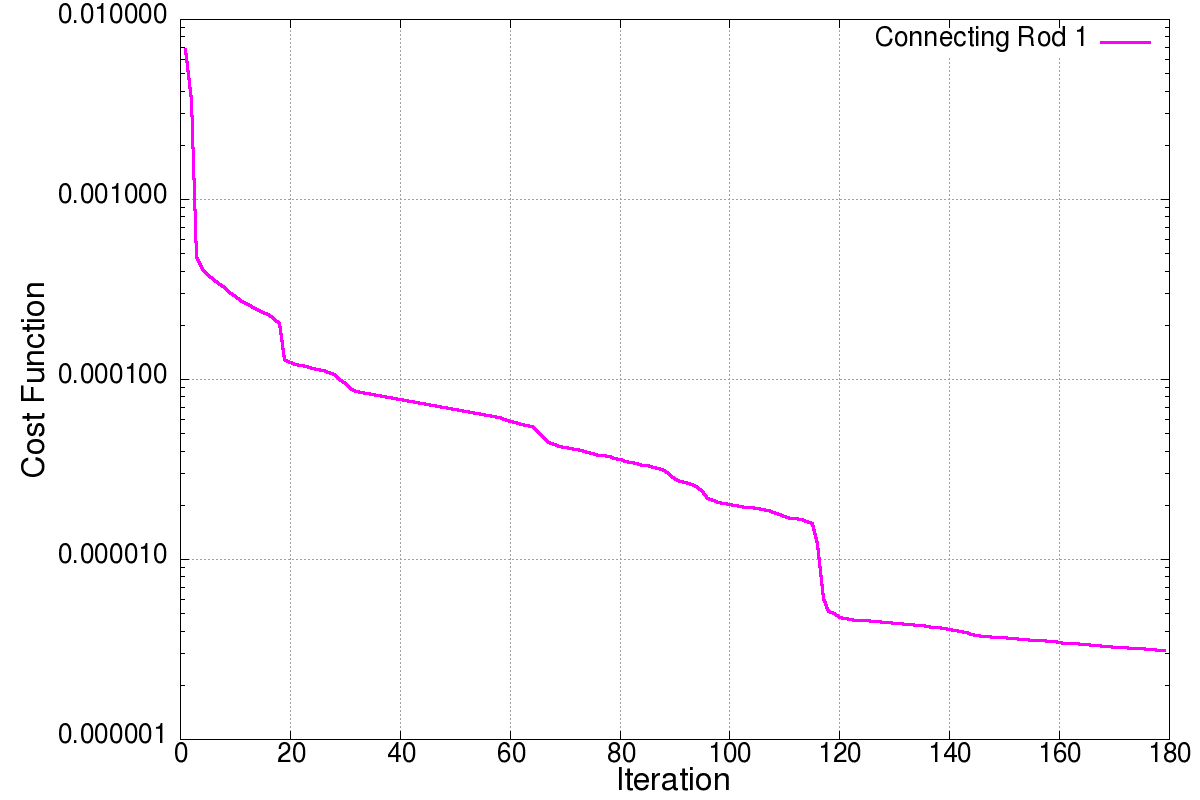}
    \caption{Connecting Rod: Cost Function History}
    \label{fig:connecting_rod_cost}
\end{figure}

\begin{figure}[!hbt]
    \centering
    \includegraphics[width=0.78\textwidth]{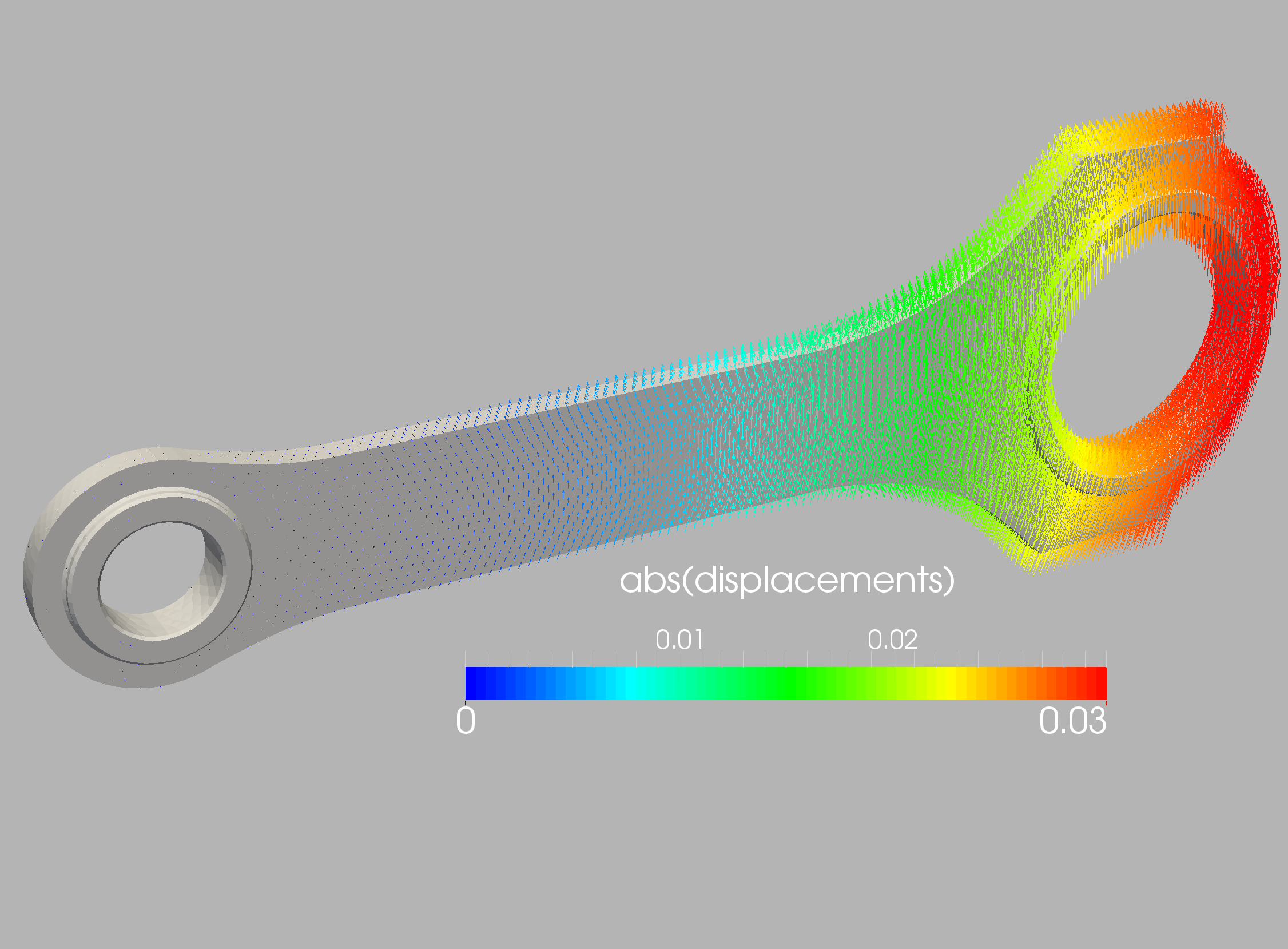}
    \caption{Connecting Rod: Displacements}
    \label{fig:connecting_rod_displacements}
\end{figure}

\begin{figure}[!hbt]
    \centering
    \includegraphics[width=0.78\textwidth]{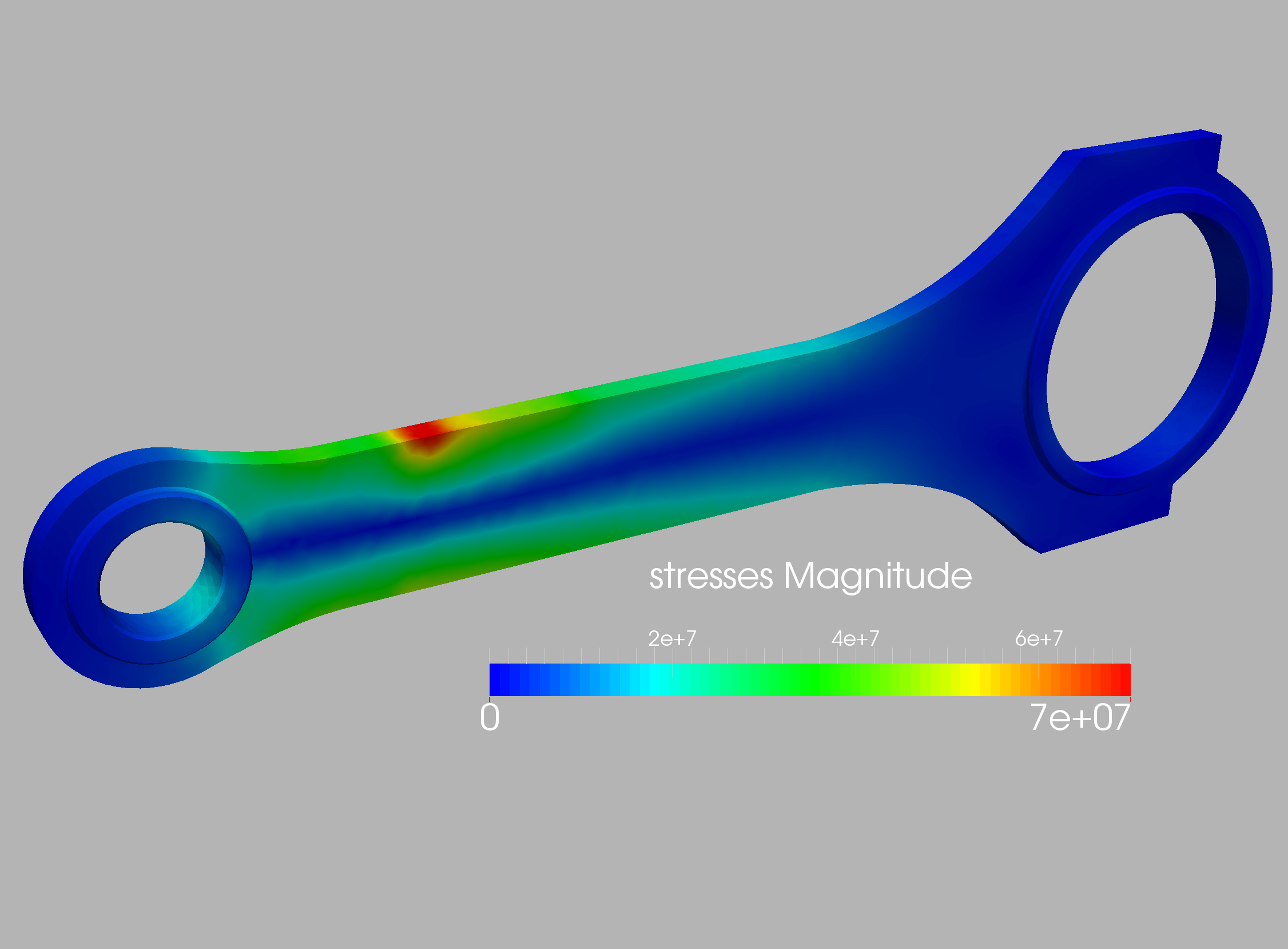}
    \caption{Connecting Rod: Stresses}
    \label{fig:connecting_rod_stresses}
\end{figure}

The case is shown in Figure \ref{fig:connecting_rod_surface} and considers a connecting rod typical
of mechanical and aerospace engineering (e.g. to actuate flaps in wings).
The two inner diameters are $d_1=2, d_2=6$, and the distance between
the centers is $d_{12}=22$ (all units in cgs).
Density, Young's modulus and Poisson rate were set to
$\rho=7.8, E=2 \cdot 10^{12}, \nu=0.3$ respectively. The inner part of
the smaller cylinder is held fixed, while forces in the $x,y$ 
direction are applied at the larger cylinder. In order to assess
the effect of mesh refinement, two different grids were employed:
coarse (9.9Ktet) and medium (71Ktet). In all cases
linear, tetrahedral elements were used. The surface mesh of the
medium mesh is shown in Figure \ref{fig:connecting_rod_wire}. In a first series of runs, 
32~measuring
points were placed on the connector rod surface and the target strength 
factor shown in Figures \ref{fig:connecting_rod_target_strf} was specified. Figures \ref{fig:connecting_rod_iter0} and \ref{fig:connecting_rod_iter160} depict the 
difference between the measured and computed displacements and the
strength factor at iterations 0 (beginning) and 160 (end).
Note the decrease in the difference between the measured and computed 
displacements, and the emergence of the weakened region. This is
also reflected in the convergence history of the cost function
(Figure \ref{fig:connecting_rod_cost}). The displacements and stresses are shown in Figures \ref{fig:connecting_rod_displacements} and \ref{fig:connecting_rod_stresses}.

\begin{figure}[!hbt]
    \centering
    \includegraphics[width=0.78\textwidth]{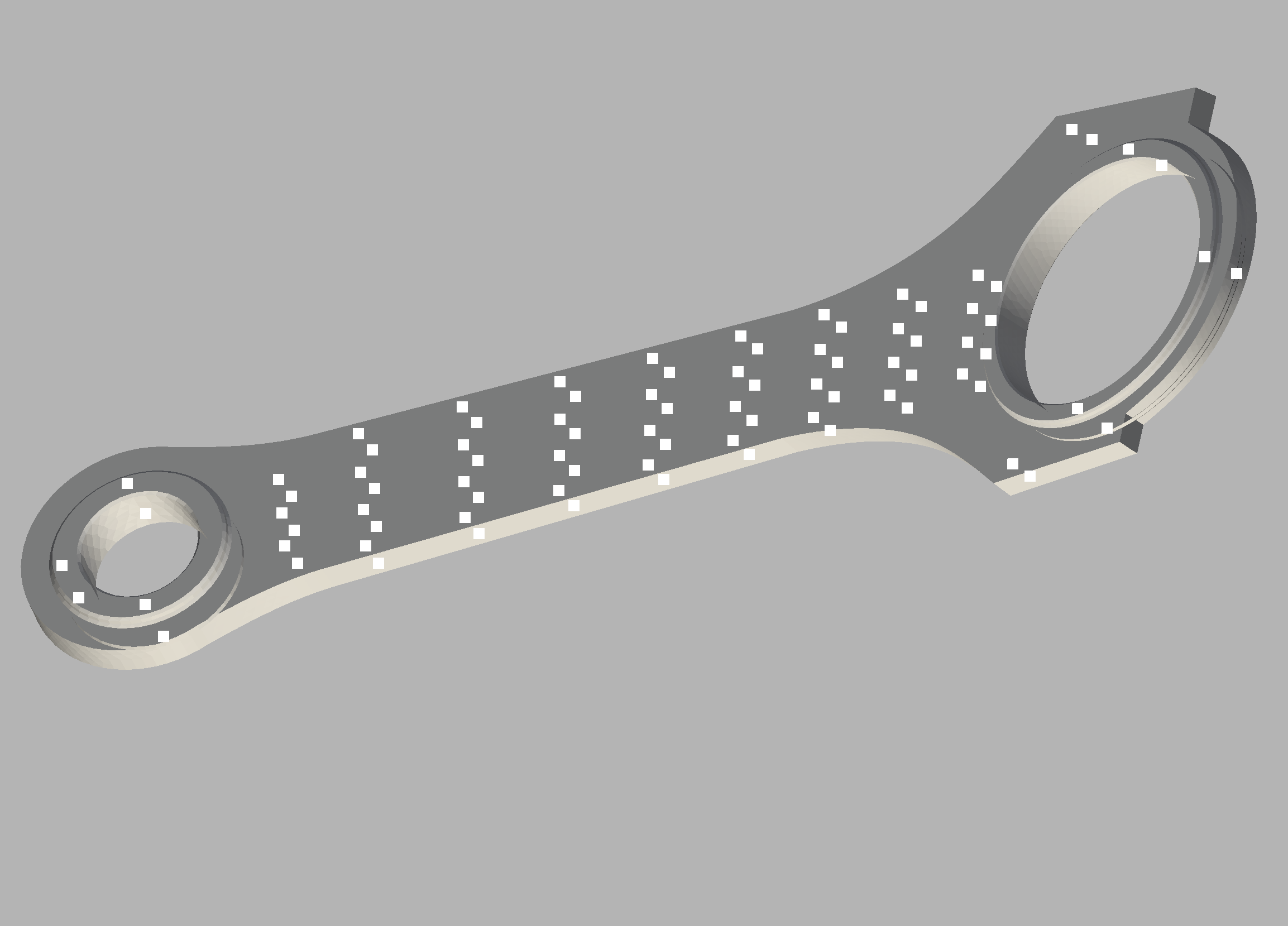}
    \caption{Connecting Rod: 86~Possible Sensor Locations}
    \label{fig:connecting_rod_possible_sensor_locations}
\end{figure}

\begin{figure}[!hbt]
    \centering
    \includegraphics[width=0.78\textwidth]{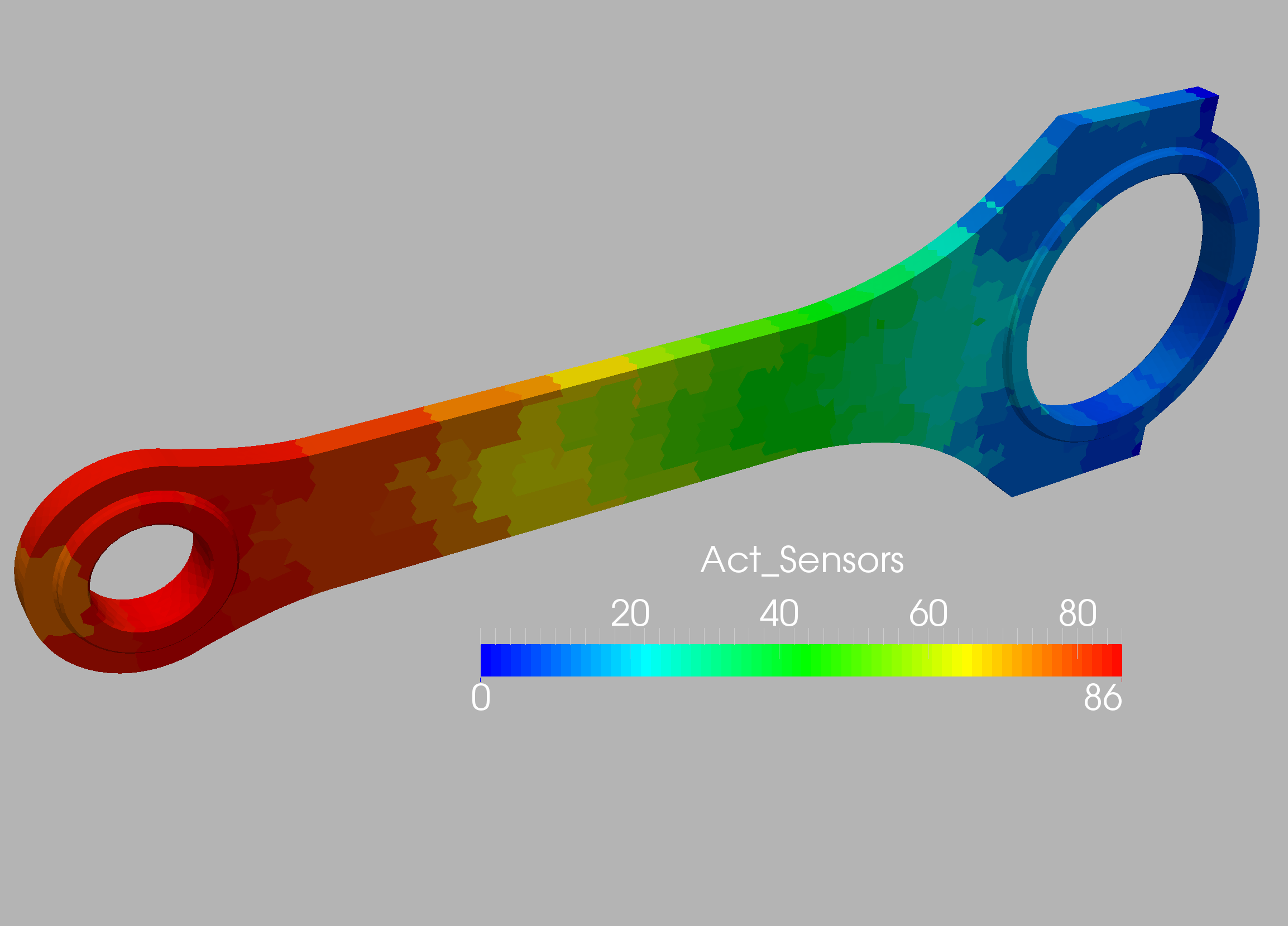}
    \caption{Connecting Rod: Number of Sensors Activated by Weakening Element Groups}
    \label{fig:connecting_rod_sensors_activated}
\end{figure}

\begin{figure}[!hbt]
    \centering
    \includegraphics[width=0.78\textwidth]{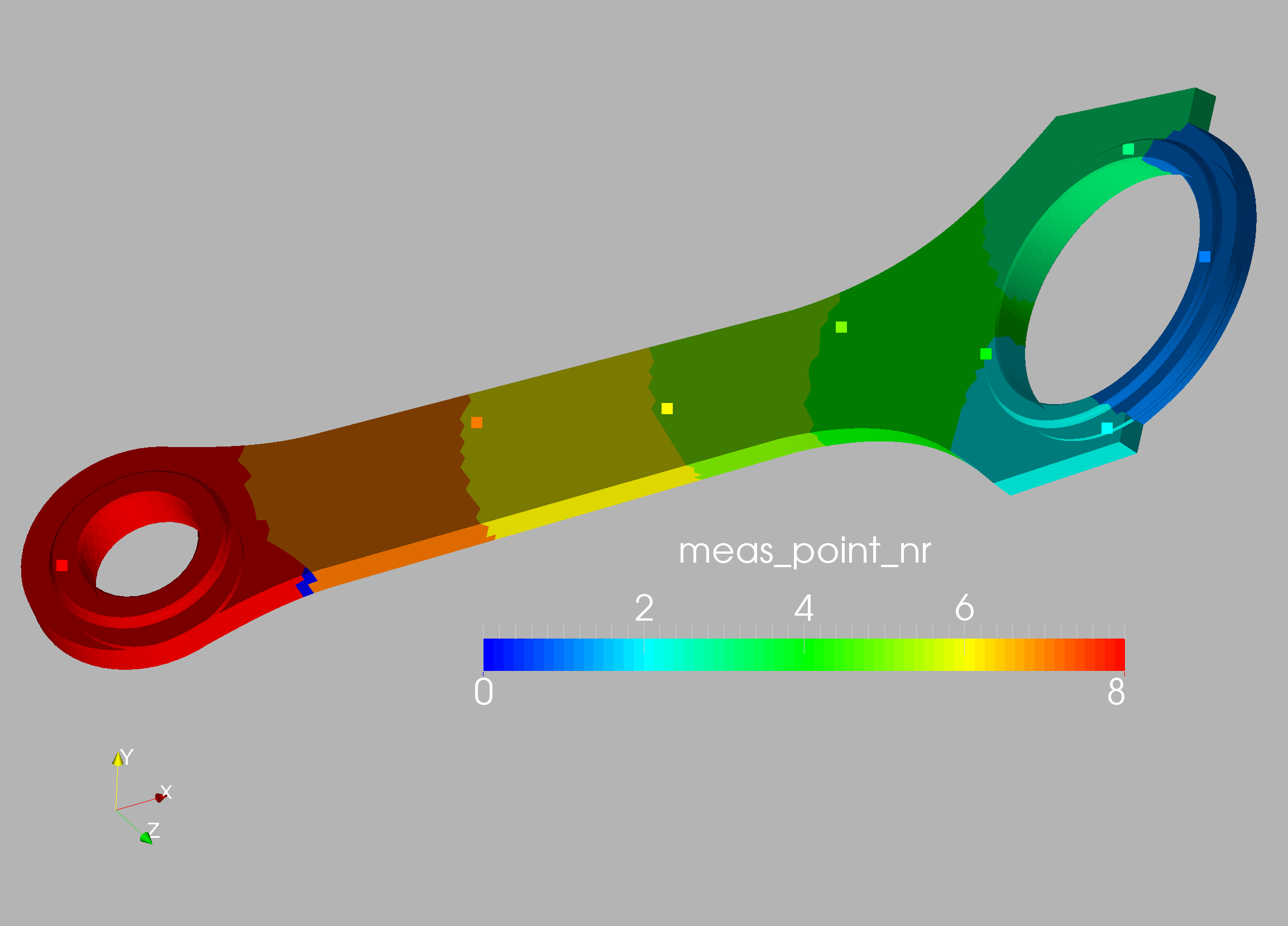}
    \caption{Connecting Rod: Optimal Sensor Locations}
    \label{fig:connecting_rod_sensor_locations}
\end{figure}

\begin{figure}[!hbt]
    \centering
    \includegraphics[width=0.78\textwidth]{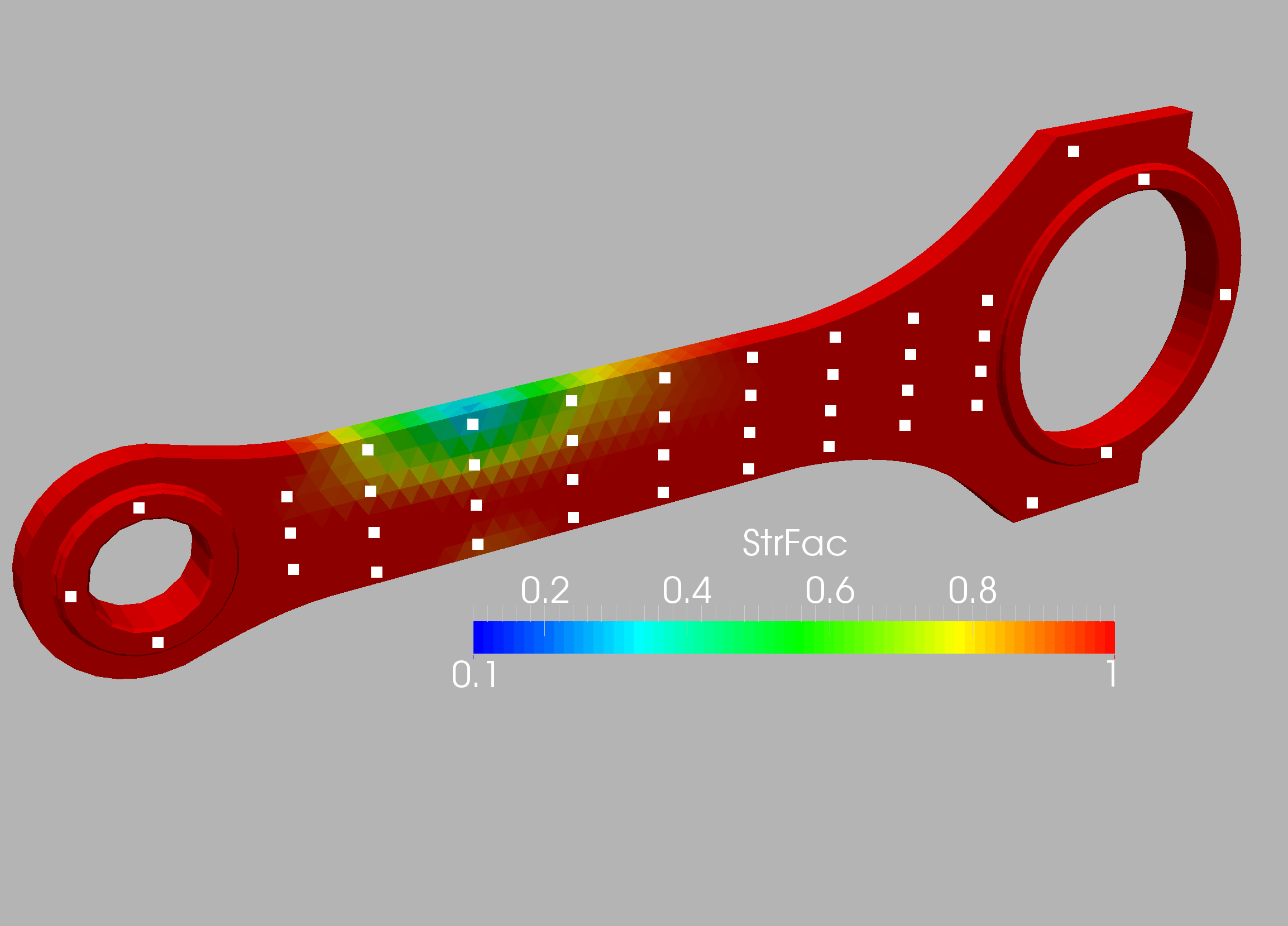}
    \caption{Connecting Rod: Strength Factor Obtained With 86 Sensors on Coarse Mesh}
    \label{fig:connecting_rod_strf_86sensor_coarse}
\end{figure}

\begin{figure}[!hbt]
    \centering
    \includegraphics[width=0.78\textwidth]{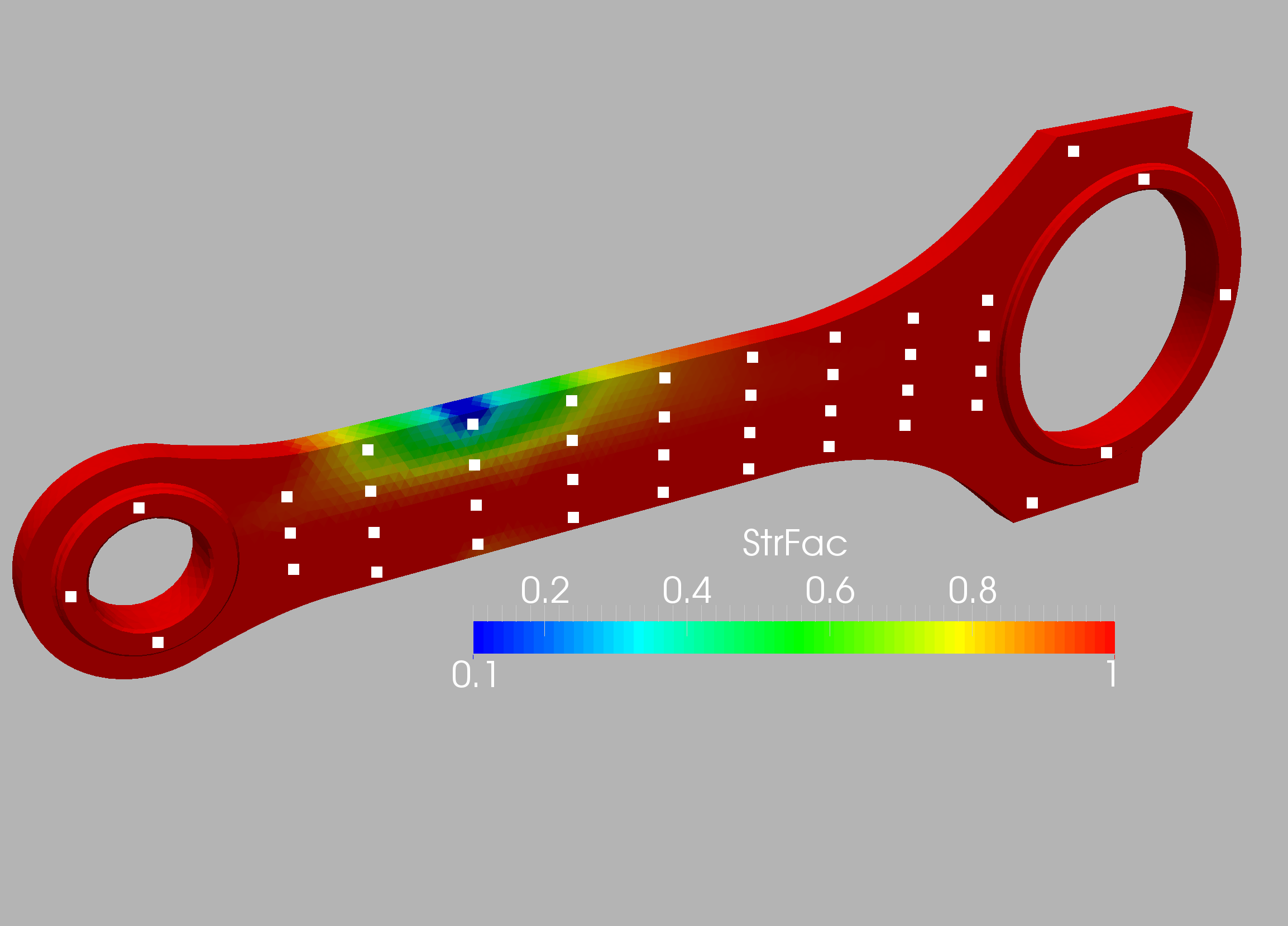}
    \caption{Connecting Rod: Strength Factor Obtained With 86 Sensors on Medium Mesh}
    \label{fig:connecting_rod_strf_86sensor_medium}
\end{figure}

\begin{figure}[!hbt]
    \centering
    \includegraphics[width=0.78\textwidth]{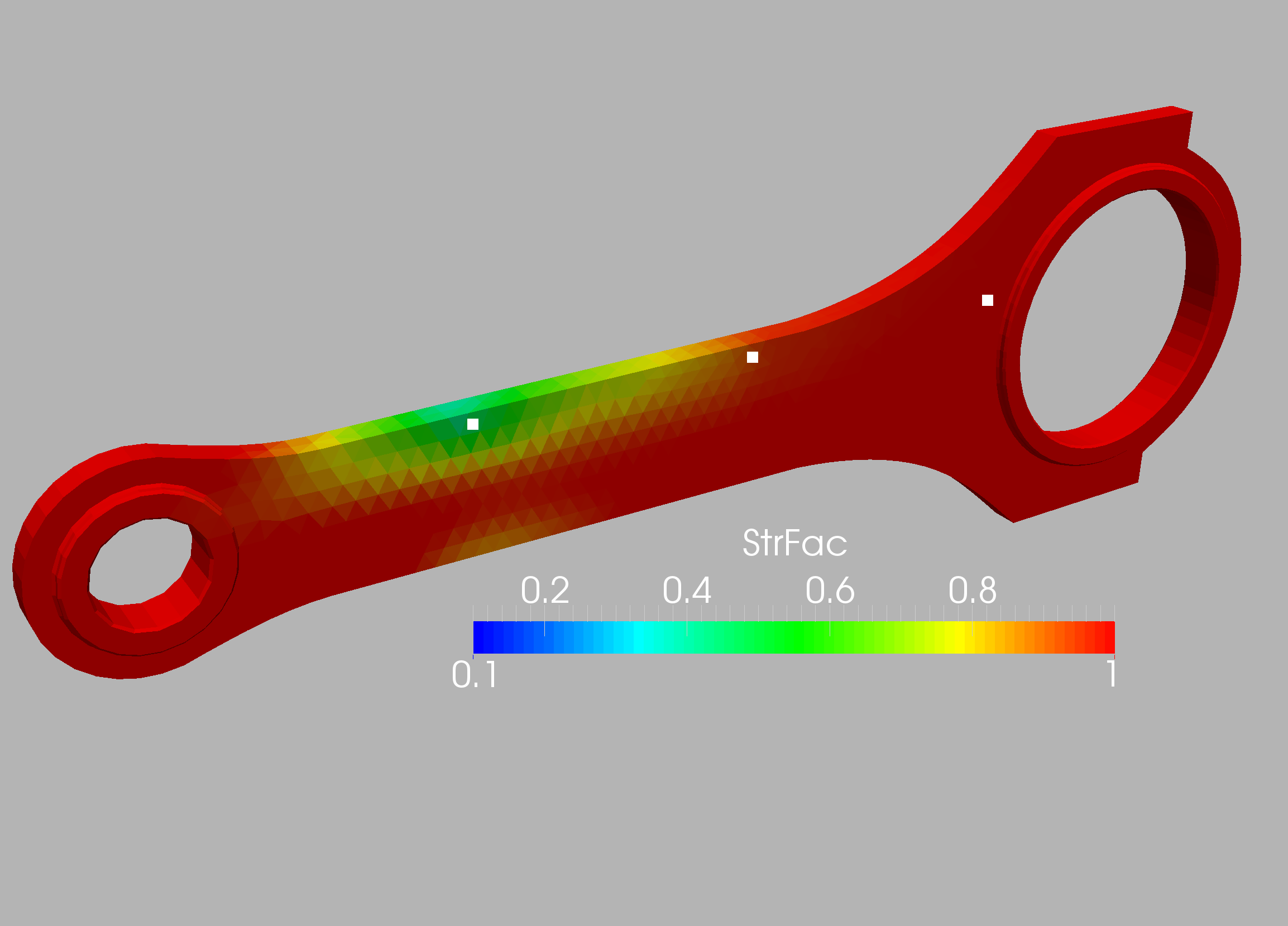}
    \caption{Connecting Rod: Strength Factor Obtained With 8 Sensors on Coarse Mesh}
    \label{fig:connecting_rod_strf_8sensor_coarse}
\end{figure}

\begin{figure}[!hbt]
    \centering
    \includegraphics[width=0.78\textwidth]{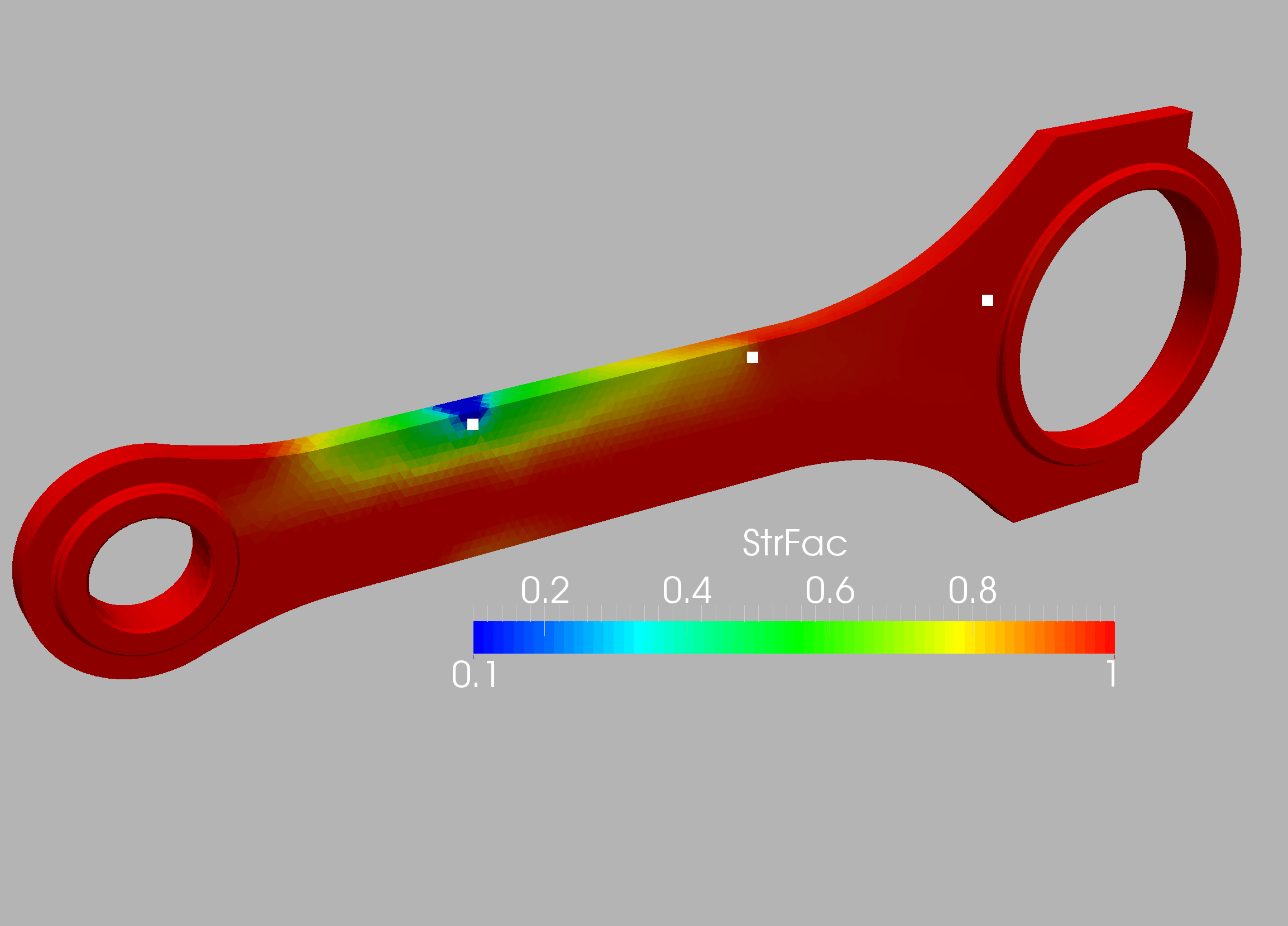}
    \caption{Connecting Rod: Strength Factor Obtained With 8 Sensors on Medium Mesh}
    \label{fig:connecting_rod_strf_8sensor_medium}
\end{figure}

\ms \noi
This case was analyzed further by specifying 86~measuring points 
(Figure \ref{fig:connecting_rod_possible_sensor_locations}) and then obtaining from these the optimal sensors using 
the procedures outlined
above. Given that the number of elements and possible sensors was 
considerable, groups of elements of desired `group size' of 1x1x1~cm
were formed. These can be discerned in Figure \ref{fig:connecting_rod_sensors_activated}, which shows the
number of sensors that would `see' (i.e. be activated) when a group
of elements is weakened. The location of the optimal sensors obtained,
as well as the region of elements each of them covers, can be seen
in Figure \ref{fig:connecting_rod_sensor_locations}. In order to assess the effect of mesh resolution the
the weakening regions obtained using the original 86~sensor locations 
and the optimal 8~sensor locations for the 
coarse and medium meshes are compared in Figures \ref{fig:connecting_rod_strf_86sensor_coarse}-\ref{fig:connecting_rod_strf_8sensor_medium}. As expected,
mesh resolution is important, and the 8~optimally placed sensors are
able to detect the weakened region with high precision. This example
highlights the importance of using high-definition digital twins and
not simpler reduced order models (ROMs) or machine learning models
(MLs) in order to localize regions of weakended material in complex
structures.

\section{Conclusions and outlook}
\label{sec:conclusions}

An adjoint-based procedure to determine weaknesses, or, more generally
the material properties of structures has been presented. Given a
series of force and displacement/strain measurements, the material
properties are obtained by minimizing the weighted differences
between the measured and computed values. In a subsequent step
techniques to optimize the number of loadings and sensors have been
proposed and tested. \\
Several examples show the viability, accuracy and efficiency of 
the proposed methodology. \\
We consider this a first step that demonstrates the viability of the
adjoint-based methodology for system identification and its use for
high fidelity digital twins \cite{mainini2015surrogate,chinesta2020virtual}. \\
Many questions remain open, of which we
just mention some obvious ones:
\begin{itemize}
\item[-] Will these techniques work for nonlinear problems ?
\item[-] Which sensor resolution is required to obtain reliable results ?
\item[-] Will these techniques work under uncertain measurements ? \cite{antil2018frontiers}
\item[-] Can one detect faulty sensors in a systematic way ?
\end{itemize}

\section{Acknowledgements}
\label{sec:acknowledgements}

This work has been partially supported by NSF grant DMS-2110263 and 
the Air Force Office of Scientific Research (AFOSR) under Award 
NO: FA9550-22-1-0248.

\bibliographystyle{plain}
\bibliography{references.bib}  

\begin{thebibliography}{10}

\bibitem{airaudo2023adjoint}
Facundo~N Airaudo, Rainald L{\"o}hner, Roland W{\"u}chner, and Harbir Antil.
\newblock Adjoint-based determination of weaknesses in structures.
\newblock {\em Computer Methods in Applied Mechanics and Engineering},
  417:116471, 2023.

\bibitem{alkayem2018structural}
Nizar~Faisal Alkayem, Maosen Cao, Yufeng Zhang, Mahmoud Bayat, and Zhongqing
  Su.
\newblock Structural damage detection using finite element model updating with
  evolutionary algorithms: a survey.
\newblock {\em Neural Computing and Applications}, 30:389--411, 2018.

\bibitem{antil2018frontiers}
Harbir Antil, Drew~P Kouri, Martin-D Lacasse, and Denis Ridzal.
\newblock {\em Frontiers in PDE-constrained Optimization}, volume 163.
\newblock Springer, 2018.

\bibitem{borrvall2003topology}
Thomas Borrvall and Joakim Petersson.
\newblock Topology optimization of fluids in stokes flow.
\newblock {\em International journal for numerical methods in fluids},
  41(1):77--107, 2003.

\bibitem{bunting2015eigenmodes}
Gregory Bunting, Scott~T Miller, Timothy~F Walsh, Clark~R Dohrmann, and Wilkins
  Aquino.
\newblock Novel strategies for modal-based structural material identification.
\newblock {\em Mechanical Systems and Signal Processing}, 149:107295, 2021.

\bibitem{bunting2021novel}
Gregory Bunting, Scott~T Miller, Timothy~F Walsh, Clark~R Dohrmann, and Wilkins
  Aquino.
\newblock Novel strategies for modal-based structural material identification.
\newblock {\em Mechanical Systems and Signal Processing}, 149:107295, 2021.

\bibitem{cawley1979location}
Peter Cawley and Robert~Darius Adams.
\newblock The location of defects in structures from measurements of natural
  frequencies.
\newblock {\em The Journal of Strain Analysis for Engineering Design},
  14(2):49--57, 1979.

\bibitem{chamoin2014updatingoffemmodels}
Ludovic Chamoin, Pierre Ladev{\`e}ze, and Julien Waeytens.
\newblock Goal-oriented updating of mechanical models using the adjoint
  framework.
\newblock {\em Computational mechanics}, 54(6):1415--1430, 2014.

\bibitem{chinesta2020virtual}
Francisco Chinesta, Elias Cueto, Emmanuelle Abisset-Chavanne, Jean~Louis Duval,
  and Fouad~El Khaldi.
\newblock Virtual, digital and hybrid twins: a new paradigm in data-based
  engineering and engineered data.
\newblock {\em Archives of computational methods in engineering}, 27:105--134,
  2020.

\bibitem{dhondt2022calculix}
Guido Dhondt.
\newblock Calculix user’s manual version 2.20.
\newblock {\em Munich, Germany}, 2022.

\bibitem{di2022data}
Daniele Di~Lorenzo, Victor Champaney, Claudia Germoso, Elias Cueto, and
  Francisco Chinesta.
\newblock Data completion, model correction and enrichment based on sparse
  identification and data assimilation.
\newblock {\em Applied Sciences}, 12(15):7458, 2022.

\bibitem{kim2004damage}
Hansang Kim and Hani Melhem.
\newblock Damage detection of structures by wavelet analysis.
\newblock {\em Engineering structures}, 26(3):347--362, 2004.

\bibitem{ladeveze1994modalfemupdate}
Pierre Ladev{\`e}ze, Djamel Nedjar, and Marie Reynier.
\newblock Updating of finite element models using vibration tests.
\newblock {\em AIAA journal}, 32(7):1485--1491, 1994.

\bibitem{lazarov2011filters}
Boyan~Stefanov Lazarov and Ole Sigmund.
\newblock Filters in topology optimization based on helmholtz-type differential
  equations.
\newblock {\em International Journal for Numerical Methods in Engineering},
  86(6):765--781, 2011.

\bibitem{lohner2008applied}
Rainald L{\"o}hner.
\newblock {\em Applied computational fluid dynamics techniques: an introduction
  based on finite element methods}.
\newblock John Wiley \& Sons, 2008.

\bibitem{lohner2023feelast}
Rainald L\"ohner.
\newblock Feelast user’s manual.
\newblock {\em Fairfax, Virginia}, 2023.

\bibitem{lohner2020determination}
Rainald L{\"o}hner and Harbir Antil.
\newblock Determination of volumetric material data from boundary measurements:
  Revisiting calderon’s problem.
\newblock {\em International Journal of Numerical Methods for Heat \& Fluid
  Flow}, 2020.

\bibitem{maia1997localization}
NMM Maia, JMM Silva, and RPC Sampaio.
\newblock Localization of damage using curvature of the
  frequency-response-functions.
\newblock In {\em Proceedings of the 15th international modal analysis
  conference}, volume 3089, page 942, 1997.

\bibitem{mainini2015surrogate}
Laura Mainini and Karen Willcox.
\newblock Surrogate modeling approach to support real-time structural
  assessment and decision making.
\newblock {\em AIAA Journal}, 53(6):1612--1626, 2015.

\bibitem{mirzaee2015adjointaccel}
Akbar Mirzaee, Reza Abbasnia, and Mohsenali Shayanfar.
\newblock A comparative study on sensitivity-based damage detection methods in
  bridges.
\newblock {\em Shock and Vibration}, 2015, 2015.

\bibitem{mohan2013structural}
SC~Mohan, Dipak~Kumar Maiti, and Damodar Maity.
\newblock Structural damage assessment using frf employing particle swarm
  optimization.
\newblock {\em Applied Mathematics and Computation}, 219(20):10387--10400,
  2013.

\bibitem{puelaubry2011meshadapt}
Guillaume Puel and Denis Aubry.
\newblock Using mesh adaption for the identification of a spatial field of
  material properties.
\newblock {\em International journal for numerical methods in engineering},
  88(3):205--227, 2011.

\bibitem{rucka2006application}
Magdalena Rucka and Krzysztof Wilde.
\newblock Application of continuous wavelet transform in vibration based damage
  detection method for beams and plates.
\newblock {\em Journal of sound and vibration}, 297(3-5):536--550, 2006.

\bibitem{salloum2022optimization}
Maher Salloum and David~B Robinson.
\newblock Optimization of flow in additively manufactured porous columns with
  graded permeability.
\newblock {\em AIChE Journal}, 68(9):e17756, 2022.

\bibitem{seidl2019forwbackw}
D~Thomas Seidl, Assad~A Oberai, and Paul~E Barbone.
\newblock The coupled adjoint-state equation in forward and inverse linear
  elasticity: Incompressible plane stress.
\newblock {\em Computer Methods in Applied Mechanics and Engineering},
  357:112588, 2019.

\bibitem{simo2006computational}
Juan~C Simo and Thomas~JR Hughes.
\newblock {\em Computational inelasticity}, volume~7.
\newblock Springer Science \& Business Media, 2006.

\bibitem{troltzsch2010optimal}
Fredi Tr{\"o}ltzsch.
\newblock {\em Optimal control of partial differential equations: theory,
  methods, and applications}, volume 112.
\newblock American Mathematical Soc., 2010.

\bibitem{zienkiewicz2005finite}
Olek~C Zienkiewicz, Robert~Leroy Taylor, and Jian~Z Zhu.
\newblock {\em The finite element method: its basis and fundamentals}.
\newblock Elsevier, 2005.

\end{thebibliography}

\end{document}